\documentclass[a4paper,12pt]{amsart}

\usepackage[latin1]{inputenc}
\usepackage[english]{babel}
\usepackage{amsmath}
\usepackage{amssymb}
\usepackage{enumerate}
\usepackage{xspace}
\usepackage{amsfonts}
\usepackage{latexsym}
\usepackage{amsthm}
\usepackage{amsrefs}
\usepackage[all]{xypic}

\newtheorem{lemma}{Lemma}

\newtheorem{theorem}[lemma]{Theorem}

\newtheorem{corollary}[lemma]{Corollary}

\theoremstyle{definition}
\newtheorem{definition}[lemma]{Definition}

\newtheorem{remark}[lemma]{Remark}

\newcommand{\HX}{\mathfrak{H}_{\OSS}}
\newcommand{\lanl}{\mathcal{L}(\Lambda)}
\newcommand{\crol}{\cro{\mathcal{D}_{\OSS_\Lambda}\mspace{-6mu}}{\alpha}{\tra}}
\newcommand{\croy}{\cro{\D_{\OSSY}\mspace{-2mu}}{\alpha}{\tra}}
\newcommand{\croa}{\cro{\D_{\OSS_A}\mspace{-4mu}}{\alpha}{\tra}}
\newcommand{\AX}{\mathcal{A}_{\OSS}}
\newcommand{\gaug}{\gamma}
\newcommand{\C}{\mathbb{C}}
\newcommand{\N}{\mathbb{N}}
\newcommand{\No}{\N_0}
\newcommand{\Z}{\mathbb{Z}}
\newcommand{\R}{\mathbb{R}}
\newcommand{\A}{\mathcal{A}}
\newcommand{\DX}{\mathcal{D}_{\OSS}}
\newcommand{\floweq}{\cong_f}
\newcommand{\lamb}[1]{\lambda_{#1}}
\newcommand{\crod}{\cro{\DX\mspace{-2mu}}{\alpha}{\tra}}
\newcommand{\quo}{\rho}
\newcommand{\cro}[3]{#1\rtimes_{#2,#3}\N}
\newcommand{\cros}{\cro{A}{\alpha}{\tra}}
\newcommand{\toep}{\mathcal{T}}
\newcommand{\toepd}{\toep(\DX,\alpha,\tra)}
\newcommand{\toepl}{\toep(A,\alpha,\tra)}
\newcommand{\tra}{\mathcal{L}}
\newcommand{\cy}[2]{C(#1,#2)}
\newcommand{\cs}{C^*}
\newcommand{\inv}{^{-1}}
\newcommand{\al}{\mathfrak a}
\newcommand{\Past}{{\mathcal P}}
\newcommand{\OSS}[1][]{{\mathsf{X}_{#1}}}
\newcommand{\osh}{\sigma}
\newcommand{\tsh}{\sigma}
\newcommand{\alwords}{\al^*}
\newcommand{\cyl}[1]{C(#1)}
\newcommand{\B}{\mathsf{L}}
\newcommand{\E}{\mathcal{E}}
\newcommand{\D}{\mathcal{D}}
\newcommand{\Oo}{\mathcal{O}}
\newcommand{\F}{\mathcal{F}}
\newcommand{\K}{\mathcal{K}}
\newcommand{\T}{\mathbb{T}}
\newcommand{\limm}{\underset{\longrightarrow}{\lim}}
\newcommand{\HH}{\mathsf{H}}
\newcommand{\emptyword}{\epsilon}
\newcommand{\OSSY}{{\mathsf{Y}}}
\DeclareMathOperator{\aut}{Aut}

\def\abstract{\vspace{0cm}
{\bf \small Abstract.}\footnotesize }

\title{$\cs$-crossed Products and Shift Spaces}

\keywords{}
\date{June 24, 2005}

\subjclass[2000]{}

\begin{document}
\maketitle

\begin{center}
Toke Meier Carlsen \\*[0.3cm]
Mathematisches Institut, 
Einsteinstraße 62, 48149 Münster, Germany \\ 
toke@math.uni-muenster.de \\*[0.7cm] 
Sergei Silvestrov\\*[0.3cm] 
Centre for Mathematical Sciences,
Department of Mathematics, \\
Lund Institute of Technology,
Lund University,\\
Box 118, 221 00 Lund, Sweden \\ sergei.silvestrov@math.lth.se \\
FAX:  +46 46 2224010 \ \ tel: +46 46 2228854 \vspace{0.5cm}
\end{center}

\begin{abstract}
In this article, we use Exel's construction to
associate a $\cs$-algebra to every shift space. We show that
it has the $\cs$-algebra defined in \cite{MR2091486} as a quotient, 
and possesses properties indicating that it can be thought of  
as the universal $\cs$-algebra associated to a shift space.
We also consider its representations, relationship to other $C^*$-algebras associated to shift spaces, show that it can be viewed as a generalization of the universal Cuntz-Krieger algebra, 
discuss uniqueness and a faithful representation, provide conditions for it 
being nuclear, for satisfying the UCT, for being simple, and for being purely infinite, 
show that the constructed algebras and thus their K-theory, $K_0$ and $K_1$, are conjugacy invariants of one-sided shift spaces, present formulas for those invariants, and  
also present a description of the structure of gauge invariant ideals.   
\\ \\   
{\bf Keywords:} $\cs$-algebra, shift spaces, dynamical systems, invariants.
\end{abstract}

\footnotetext[1]{Mathematics Subject
Classification 2000: Primary 47L65; Secondary 46L55, 37B10, 54H20}

\footnotetext[2]{Supported by The Swedish Foundation for
International Cooperation in Research and High Education STINT and by the Crafoord Foundation.}

\section{Introduction}
When dynamical system consists of a homeomorphism of a  topological space, or more generally when an action of a group of invertible transformations of some space is studied, there is a standard
construction of a crossed product $\cs$-algebra. Historically this
construction has its origins in foundations of quantum mechanics.
The important idea behind this construction is that it encodes the
action and the space within one algebra thus providing
opportunities for their investigation on the same level. It is
known that properties of the topological space can be considered
via properties of the algebra of continuous functions defined on
it. The crossed product algebra is constructed by combining this algebra of functions with the action being encoded using further
elements of the new in general non-commutative algebra. The action
is built into multiplication in the new algebra via covariance
commutation relations between the elements in  the algebra of
functions and the elements used to encode the action. The crossed
product construction have considerable applications in quantum
mechanics and quantum field theory, 
and provide an important source
of examples for further development of non-commutative geometry. A
lot of research has been done on interplay between properties of
the invertible dynamical systems and properties of the
corresponding crossed product $\cs$-algebras and $W^*$-algebras.

There are several ways to generalize the construction of the
$\cs$-crossed product to the non-invertible setting. The one we will
focus on in this paper was introduced by Exel in
\cite{MR2032486}. This construction relies on a choice of
\emph{transfer operator}. Exel showed that for a natural choice of
transfer operator, the $\cs$-algebra of a one-sided shift of finite
type is isomorphic a \emph{Cuntz-Krieger algebra}.

The Cuntz-Krieger algebras was introduced by Cuntz and Krieger in
\cite{MR82f:46073a}. They can in a natural way be viewed as universal
$C^*$-algebras associated with shift spaces (also called subshifts) of finite
type. From the point of view of
operator algebra these $C^*$-algebras were important examples of
$C^*$-algebras with new properties and from the point of view of
topological dynamics these $C^*$-algebras (or rather, the $K$-theory
of these $C^*$-algebras) gave new invariants of shift spaces of
finite type.

In \cite{MR98h:46077} Matsumoto tried to generalize this idea by
constructing $C^*$-algebras associated with every shift space and
he studied them in
\cites{MR2000e:46087,MR2000d:46082,MR2000f:46084,MR2001e:46115,MR2001g:46147}.
Unfortunately there is a mistake in \cite{MR2000f:46084} which
makes many of the results in
\cites{MR2000e:46087,MR2000d:46082,MR2000f:46084,MR2001e:46115,MR2001g:46147}
invalid for the $\cs$-algebra constructed in \cite{MR98h:46077},
and since this mistake was discovered, there has been some
confusion about the right definition of the $\cs$-algebra
associated to a shift space.

In this paper we will use Exel's construction to
associate a $\cs$-algebra to every shift space, and we will show that
it has the properties Matsumoto thought his algebra had, and thus that
it satisfies all the results of
\cites{MR98h:46077,MR2000e:46087,MR2000d:46082,MR2000f:46084,MR2001e:46115,MR2001g:46147}
and has the $\cs$-algebra defined in \cite{MR2091486} as a quotient.
Thus it seems right to think of this $\cs$-algebra as the universal
$\cs$-algebra associated to a shift space.

Matsumoto's original construction associated a $\cs$-algebra to
every \emph{two-sided} shift space, but it seems more natural to
work with \emph{one-sided} shift spaces, so we will do that in
this paper, but since every two-sided shift space comes with a
canonical one-sided shift space (see below), the $\cs$-algebras we
define in this paper can in a natural way also be seen as
$\cs$-algebras associated to two-sided shift spaces.


\section{$\cs$-algebras of invertible dynamical systems}
In this section we review the construction and some properties of
a $C^*$-crossed product of a $C^*$-algebra by the action of the
discrete group of automorphisms. In particular the invertible
dynamical systems generated by homeomorphisms of topological
spaces are encoded in the crossed product $C^*$-algebras obtained
from the actions of the group of integers on the $C^*$-algebra of complex-valued continuous functions.

Let $(A,G,\alpha)$ be a triple consisting of a unital
$\cs$-algebra, discrete group $G$ and an action $\alpha: G
\rightarrow \aut(A)$ of $G$ on $A$, meaning a homomorphism from
the group $G$ into the group $\aut(A)$ of automorphisms of the $\cs$-algebra $A$. A pair $\{\pi, u\}$ consisting of a representation $\pi$ of $A$ and a unitary representation $u$ of $G$ on a Hilbert space $H$ is called a
covariant representation of the system $(A,G,\alpha)$ if
$$u_s \pi(a) u_s^*= \pi (\alpha_s(a))$$ for every $a\in A$ and $s\in G$. The full crossed product ${A}\rtimes_{\alpha} G$ is defined as the universal $C^*$-algebra for the family of covariant representations.

Another more concrete way to define ${A}\rtimes_{\alpha} G$ is to
consider the space $l^1(G,A)$ of all $A$-valued functions
$x(\cdot)$ on $G$ with the finite $l_1$-norm $||x||=\sum_{s\in G}
||x(s)||_{A}$ equipped with the twisted convolution product and
the involution
$$xy(s) = \sum_{t\in G} x(t) \alpha_t(y(t^{-1}s)), \quad
x^*(s)=\alpha_s(x(s^{-1})^*)$$ making $l^1(G,A)$ into a Banach
$*$-algebra. The algebra $A$ can be identified with the algebra of
functions $\tilde{a}: G \rightarrow A$ defined as $\tilde{a}(e) =
a \in A$ on the unit element $e$ of $ G $ and as zero elsewhere
on $G$. Moreover, for each $s\in G$ a function $\delta_s: G
\rightarrow A$ is defined as zero everywhere on $G$ except $s$
where $\delta_s(s)=1_A$ the unit element of $A$. With this
notation $\tilde{a} =a \delta_e$. It can be shown that the
functions $\delta_s$, $s\in G$ are unitary elements of the Banach
$*$-algebra $l^1(G,A)$, that is
$\delta_s\delta_s^*=\delta_s^*\delta_s = 1_{l^1(G,A)} = \delta_e$;
the map $s \mapsto \delta_s $ is a group homomorphism in the sense
that $\delta_{uv}=\delta_{u}\delta_{v}$; and moreover the
covariance relation
$$\delta_s \tilde{a} \delta_s^* = \alpha_s (\tilde{a})$$
holds for every $a\in A $ and $s\in G$. When the functions $x \in
l^1(G,A)$ are expressed as $x=\sum_{s\in G} x(s) \delta_s$ the
covariance relation implies that the operations of twisted product
and involution in $l^1(G,A)$ are the natural ones.

It can be shown that the Banach $*$-algebra $l^1(G,A)$  has
sufficiently many representations (i. e. for any $a\in l^1(G,A)$
there is a representation $\pi$ with $\pi(a)\neq 0$). Thus one can
define the $C^*$-envelope $C^*(l^1(G,A))$ as the completion of
$l^1(G,A)$ with the norm 
$$||x||_{\infty} = \sup \{||\tilde{\pi}
(x) ||\mid \tilde{\pi} \text{ is representation of } l^1(G,A)\}.$$
Any covariant representation $\{\pi, u\}$ yields a representation
$\tilde{\pi}$ of $l^1(G,A)$, and hence of $C^*(l^1(G,A))$, defined
by
$$\tilde{\pi}(x)= \sum_{s\in G} \pi(x(s))u_s$$
for $x$ with finite support (i. e. zero outside a finite subset of
$G$). Moreover, any representation of $C^*(l^1(G,A))$ has the
above form. So, $C^*(l^1(G,A))$ is the same as the full
$C^*$-crossed product ${A}\rtimes_{\alpha} G$. It is also useful
to have in mind that the subspace of finite sums $\{\sum_{s\in J}
a_s \delta_s \mid J \text{ is finite, } a_s \in A \}$ is a dense
$*$-subalgebra of ${A}\rtimes_{\alpha} G$.

Suppose that $A$ is acting on a Hilbert space $H$ and write the
action as $a h$ for $a\in A$ and $h \in H$. Let $K = l^2(G)\otimes
H$ be regarded as $l^2(G,H)$, the space of $H$-valued
$l^2$-functions on $G$ with values in $H$. A pair
$\{\pi_{\alpha},\lambda\}$ consisting of the representation
$\pi_{\alpha}$ of $A$ and a unitary representation $\lambda:
s\mapsto \lambda_s$ of $G$ on $K$ defined by
\begin{eqnarray*}
(\pi_{\alpha}(a) f)(s) &=& \alpha_{s^{-1}}(a) f(s), \quad f\in K, a \in A \\
(\lambda_s f)(t)       &=& f(s^{-1} t)
\end{eqnarray*}
is a covariant representation. The reduced crossed product
${A}\rtimes_{\alpha,r} G$ is the $C^*$-algebra acting on $K$
generated by the operator family $\{\pi_{\alpha}(a), \lambda_s
\mid a \in A, s \in G \}$. It can be proved that the definition
does not depend on the space $H$. The reduced and full crossed
products are isomorphic if  and only if the group $G$ belongs to a
class of so called amenable groups. In particular the group
$G=\mathbb{Z}$ of special relevance in connection to invertible
dynamical systems belongs to this class.

When $G=\mathbb{Z}$, the number $1\in \mathbb{Z}$ is the generator
of the group $\mathbb{Z}$. As $s\mapsto \alpha_s$ is a
homomorphism, it is enough to specify the defining covariance
relation for ${A}\rtimes_{\alpha} G$ for the generator of
$\mathbb{Z}$, that is
$$ \delta_1 a \delta_1^*=\alpha_1(a).$$

An object of special interest to us is the crossed product
$C^*$-algebra for an invertible dynamical system consisting of
iterations of a homeomorphism acting on a topological space.

Let $\Sigma = (X,\sigma)$ be a topological dynamical system
consisting of a homeomorphism of a Hausdorff topological space
$X$. The $*$-algebra of all continuous functions on $X$ and the
$*$-algebra of all continuous functions on $X$ with compact
support will be denoted respectively by $C(X)$  and by $C_{c}(X)$.
The algebra $C(X)$ has a unit if and only if $X$ is compact, and
the unit then is the constant function ${\bf 1} ={\bf
1}_{C(X)}(\cdot)$ equal to $1$ on all elements of $X$. Moreover,
$X$ is compact if and only if $C(X)$ and  $C_{c}(X)$ coincide.

The mapping $\alpha : C(X) \rightarrow C(X) $ defined by
\begin{equation} \label{eq:aut}
\alpha(f)(x) = f(\sigma^{-1}(x))
\end{equation}
is an automorphism of the $*$-algebra $C(X)$, and the mapping
defined by
\begin{equation} \label{eq:homaut}
j \mapsto \alpha^j(f)(x) = f(\sigma^{-j}(x))
\end{equation}
is a homomorphism of $\mathbb{Z}$ into the group $\aut(C(X))$ of $*$-automorphisms of $C(X)$. 
Since $\sigma$ is a homeomorphism, the
family of all compact subsets of $X$ is invariant with respect to
$\sigma$ and $\sigma^{-1}$, and  hence $\alpha$ leaves 
the $*$-subalgebra $C_{c}(X)$ of $C(X)$ invariant. The group $\mathbb{Z}$ is
a locally compact group with respect to the discrete topology,
i.e. the topology where any subset of $\mathbb{Z}$ is open. A
subset of $\mathbb{Z}$ is compact if and only if it is finite. The
set $C_{c}(\mathbb{Z},C(X))$ of continuous mappings from
$\mathbb{Z}$ to $C(X)$ with compact support consists of all
mappings which may assume non-zero values only at finitely many
elements of $\mathbb{Z}$. For any function $a : \mathbb{Z}
\rightarrow C(X) $ we denote by $a[k]$ the element of $C(X)$ equal
to the value of $a$ at $k \in \mathbb{Z} $. The pointwise addition
and multiplication by complex numbers makes
$C_{c}(\mathbb{Z},C(X))$ into a linear space, which becomes a
normed $*$-algebra with the multiplication, involution and norm
defined by
\begin{eqnarray}
(ab)[k](\cdot) & = & \sum_{s\in \mathbb{Z}} a[s](\cdot)
   \alpha^s(b[-s+k])(\cdot) =  \label{mult:crossprod} \\
  & = & \sum_{s\in \mathbb{Z}} a[s](\cdot)b[-s+k](\sigma^{-s}(\cdot)), 
  \nonumber \\
b^*[k](\cdot) &=& \alpha^k (\overline{b[-k]})(\cdot)
              = \overline{b[-k]}(\sigma^{-k}(\cdot)), 
              \label{inv:crossprod} \\
\| b \| & = & \sum_{s\in \mathbb{Z}} \| b[s] \|_{C(X)}.
\label{norm:crossprod}
\end{eqnarray}
The Banach $*$-algebra obtained as the completion of this normed
$*$-algebra is denoted by $l^{1}(\mathbb{Z},C(X))$.

Let us assume that $X$ is compact. Then $C_{c}(X)$ coincides with
$C(X)$. The $*$-algebra $C(X)$ becomes a unital $C^*$-algebra with
respect to the supremum norm defined by $ \| f \| = \| f \|_{C(X)}
= \sup \{f(x) \mid x \in X \} $ for all $f \in C(X)$. The mappings
defined by
$$\delta_j [k](\cdot)=\left\{
\begin{array}{lll}
{\bf 1} ={\bf 1}_{C(X)}(\cdot) & {\textstyle if } &  k=j \\
 0                             & {\textstyle if } &  k \neq j
\end{array}  \right. $$
for $j\in \mathbb{Z}$ belong to $C_{c}(\mathbb{Z},C(X))$, and
$\delta_{0}$ is the unit of $C_{c}(\mathbb{Z} ,C(X))$ and hence of
$l^{1}(\mathbb{Z},C(X))$. With the multiplication defined by
\mbox{(\ref{mult:crossprod})}, the equality $ \delta_{j} =
\delta_{1}^{j} $ holds for all $j \in \mathbb{Z}\setminus \{0\} $.
In what follows, for the brevity of notations, we will denote
$\delta_{1}$ by $\delta$, will assume that $\delta^{0} =
\delta_{0} $, and will write $\delta^{j}$ instead of $\delta_{j}$
for all $j \in \mathbb{Z} $. The algebra $C_{c}(\mathbb{Z},C(X))$
then coincides with the algebra of polynomials in $ \delta $ with
coefficients in $C(X)$.

The $C^*$-algebra $C(X)$ can be shown to be isomorphic to the
$C^*$-algebra $C(X) \delta^{0} $ inside the normed $*$-algebras
$C_{c}(\mathbb{Z}, C(X))$ and $l^{1}(\mathbb{Z},C(X))$ having the
same unit $\delta^{0}$. The mapping $i_{0}: C(X) \rightarrow C(X)
\delta^{0}$ sending $ f \in C(X)$ to $ f \delta^{0} \in
C_{c}(\mathbb{Z}, C(X))$ is a unital $*$-isomorphism of
the $C^*$-algebra $C(X)$ onto the $C^*$-algebra $C(X) \delta^{0} $. We use the 
notation
$$ (f\delta^{0})[k](x) = (\delta^{0}f)[k](x)=
\left\{
\begin{array}{ll}
f(x), & k=0 \\
0, & k \neq 0
\end{array}
\right. .$$ In general, whenever it is convenient, for $a \in
l^{1}(\mathbb{Z},C(X))$ and $f \in C(X)$, by equalities of the
form $a=f$ we will mean $a=f\delta^{0}$, and the  notations $af =
a(f\delta^{0})$ and $fa = (f\delta^{0})a$ will be used with
products between $a$ and $f\delta^{0}$ defined by
\mbox{(\ref{mult:crossprod})}. The same notations will often be used for $a$ belonging to the $C^*$-crossed product algebra of
$C(X)$ by $\mathbb{Z}$ obtained as the completion of
$l^{1}(\mathbb{Z},C(X))$ with respect to a certain norm. With
this notation, the fundamental equality
\begin{equation} \label{covrelalg}
\delta f \delta^{*} = \alpha (f),
\end{equation}
called the {\it covariance relation}, holds for all $f\in C(X)$.

The mapping $E : l^{1}(\mathbb{Z}, C(X)) \rightarrow C(X)
\delta^{0}$ defined by $E(b)=b[0] \delta^{0}$ for any element $b
\in l^{1} (\mathbb{Z}, C(X))$ is a projection of norm one
satisfying
\begin{eqnarray}
& E(a b c)  =  a E(b)c \mbox{ for all } a,c\in C(X) \delta^{0}, &
\quad  \mbox{ ({\em module property}) } \\
& E(b^{*}b)  \geq  0, &
\quad   \mbox{ ({\em positivity}) } \\
& E(b^{*}b) = 0 \mbox{ implies that } b=0 & \quad   \mbox{ ({\em
faithfulness}) }
\end{eqnarray}
for all $ b \in l^{1}(\mathbb{Z}, C(X)) $. The  positivity, for
example, is proved as follows:
\begin{eqnarray*}
E(b^{*}b) & = & (b^{*}b)[0] \delta^{0} =
(\sum_{k\in\mathbb{Z}}b^{*}[k](\cdot)
\alpha^{k}(b[-k])(\cdot))\delta^{0} \\
&=& (\sum_{k\in\mathbb{Z}}
\alpha^{k}(\overline{b[-k]}b[-k])(\cdot))\delta^{0}
= (\sum_{k\in\mathbb{Z}}\alpha^{k}(|b[-k]|^{2}(\cdot)))\delta^{0}\\
&=& \sum_{k\in\mathbb{Z}}
(|b[-k](\sigma^{-k}(\cdot))|^{2})\delta^{0} \geq 0
\end{eqnarray*}
where the sums converge in norm.

For any linear functional $\varphi $ on $C(X)$, the mapping
$\varphi \circ i_{0}^{-1} $ is a linear functional on $C(X)
\delta^{0} $ satisfying $ (\varphi \circ i_{0}^{-1}) (i_{0} (a)) =
\varphi (a) $ for any $a \in C(X) $. Since the mapping $a \mapsto
i_{0} (a) $ is an isometric $*$-isomorphism of $C(X)$ onto $C(X)
\delta^{0}$, it follows that  $\| \varphi \circ i_{0}^{-1} \| = \|
\varphi \|$ for any bounded $\varphi $ on $ C(X)$, and that
$\varphi $ is positive on $C(X)$ if and only if  $\varphi \circ
i_{0}^{-1} $ is positive on $C(X) \delta^{0} $.

For any positive linear functional $\varphi$ on $C(X)$, the
mapping $(\varphi \circ i_{0}^{-1})\circ E$  is a positive linear
functional on $l^{1}(\mathbb{Z}, C(X))$. Moreover,
$\|\varphi\|=\varphi(e)$ for any positive linear functional
$\varphi$ on a Banach $^{*}$-algebra with the unit $e$. Since
\begin{eqnarray*}
& \| (\varphi \circ i_{0}^{-1})\circ E \| = (\varphi \circ
i_{0}^{-1})(E (\delta^{0})) =
(\varphi \circ i_{0}^{-1})(\delta^{0}) = \\
& = (\varphi \circ i_{0}^{-1})(i_{0} ( {\bf 1}_{C(X)} )) = \varphi
( {\bf 1}_{C(X)}) = \| \varphi \| ,
\end{eqnarray*}
the functional $\varphi$ is a state on $C(X)$, i.e. a positive
linear functional with $\|\varphi\|=1$, if and only if $(\varphi
\circ i_{0}^{-1}) \circ E$ is a state on $l^{1}(\mathbb{Z},
C(X))$.

A set of states on a Banach $*$-algebra $A$ is said to contain
sufficiently many states if for any non-zero $a\in A$ there exists
a state $\varphi$ from this set such that $\varphi(a^{*}a)\neq 0$.

There are sufficiently many states on any $C^*$-algebra, and in
particular on $C(X)$ and on its isomorphic copy $C(X) \delta^{0}$.
By faithfulness of the projection $E$, the set
$$\{(\varphi \circ i_{0}^{-1})\circ E \mid \varphi
\text{ is a state on } C(X)\}$$ of states on $l^{1}(\mathbb{Z},
C(X))$ contains sufficiently many states. As a Banach $*$-algebra
with sufficiently many states, $l^{1}(\mathbb{Z}, C(X))$ has
sufficiently many representations, i.e. for any non-zero $b\in
l^{1}(\mathbb{Z}, C(X))$ there is a representation $\pi$ with
$\pi(b)\neq 0$. Then it can be shown that
$$\|b\|_{\infty} = \sup\{\|\pi(b)\| \mid \pi
\text{ is a representation of } l^{1}(\mathbb{Z}, C(X)) \} \leq \|
b \|_{1}$$ defines a $C^{*}$-norm on $l^{1}(\mathbb{Z}, C(X))$.

The $C^*$-algebra obtained as the completion of
 $l^{1}(\mathbb{Z}, C(X))$ with respect to
the norm $\| \cdot \|_{\infty}$ is called the $C^{*}$-crossed
product of $C(X)$ by $\mathbb{Z}$ with respect to the action of
$\alpha$, or the transformation group $C^{*}$-algebra associated
with the dynamical system $\sum=(X,\sigma)$. Depending on which of
those two terminologies used, this algebra is denoted either by
$C(X) \rtimes_{\alpha} \mathbb{Z}$ or by $A(\sum)$. The
$C^{*}$-algebra $A(\sum)$ coincides with the closed linear span of
all polynomial expressions built of $\delta $,
$\delta^{*}=\delta^{-1} $ and also of elements from $C(X)$, or to
be more precise from $C(X)\delta^{0}$. Because of the covariance
relation (\ref{covrelalg}), all $\delta $ and
$\delta^{*}=\delta^{-1} $ in any such polynomial expression can be
moved to the right of all elements of $C(X)$. Thus any polynomial
expression built of $\delta $, $\delta^{*}=\delta^{-1} $ and of
elements of $C(X)$ is equal to a generalized polynomial in $\delta
$, that is to an element of the form $ \sum_{j=-n}^{j=n}
f_{j}\delta^{j} $. Consequently, the $C^{*}$-algebra $A(\sum)$ can
be viewed as a closed linear span of generalized polynomials in
$\delta$ over $C(X)$. The projection $E$ can be extended from
$l^{1}(\mathbb{Z}, C(X))$ to $A(\sum)=C(X) \rtimes_{\alpha}
\mathbb{Z}$ with the property of being faithful and with $\| E \|
=1$. For an element $a$ of $A(\sum)$, the \mbox{$n$'th}
generalized Fourier coefficient $a(n)$ is defined as
$E(a(\delta^{*})^n)$.

If $\pi$ is a $*$-representation of the $C^*$-algebra $A(\Sigma)$ on a
Hilbert space $ H_{\pi}$, then $\pi^{\prime} = \pi \circ i_{0}$ is
a $*$-representation of the $C^*$-algebra $C(X)$ on $ H_{\pi}$.
If one is given  a $*$-representation $\pi^{\prime}$ of the $C^*$-algebra $C(X)$ on $ H_{\pi}$,
then $\pi = \pi^{\prime} \circ i_{0}^{-1}$ 
is a $*$-representation of
$C(X) \delta^{0}$ on $ H_{\pi}$. Moreover, $\pi^{\prime} (f) =
(\pi \circ i_{0}) (f) = \pi (f \delta^{0}) $ for any $f \in C(X)$.
With this in mind, for simplicity  of notations, if $\pi$ is a
$*$-representation of the $C^*$-algebra $A(\Sigma)$ and $f\in C(X)$, then by $\pi(f)$ we
will always mean $\pi (f \circ \delta^{0})$.

If $\pi$ is a $*$-representation  of the $C^{*}$-algebra $A(\sum)$ on
a Hilbert space $H_{\pi}$, then the unitary operator
$u=\pi(\delta)$ and the commutative set (algebra) of bounded
operators $\pi(C(X))$ on $H_{\pi}$ satisfy the set of commutation
relations
\begin{equation} \label{covrelrep}
u \pi(f) u^{*} = \pi (\alpha(f))
\end{equation}
called covariance relations or covariance relations for a set of
operators, as they are obtained by applying the $*$-representation
$\pi$ to both sides of the covariance relation (\ref{covrelalg})
in the algebra $A(\sum)$. A pair $(\pi,u)$ consisting of a
$*$-representation of the $C^{*}$-algebra $C(X)$ on a Hilbert space
$H$, and a unitary operator $u$ on $H$ satisfying the covariance
relations (\ref{covrelrep}) is called a covariant $*$-representation or simply a covariant representation 
of the system $(C(X),\alpha , \mathbb{Z})$. So, any $*$-representation
of the $C^{*}$-algebra $A(\sum)=C(X) \rtimes_{\alpha} \mathbb{Z}$
gives rise, via restriction, to a covariant 
representation of the
system $(C(X),\alpha , \mathbb{Z})$. Moreover, this covariant
representation of $(C(X),\alpha , \mathbb{Z})$ defines uniquely
the $*$-representation of the $C^{*}$-algebra $A(\sum)$, and every covariant
representation of the system $(C(X),\alpha , \mathbb{Z})$ is
obtained by restriction from a $*$-representation of the
$C^{*}$-algebra $A(\sum)=C(X) \rtimes_{\alpha} \mathbb{Z}$. In
other words, there is a one-to-one correspondence between covariant representations of the system 
$(C(X),\alpha , \mathbb{Z})$, and
$*$-representations of the $C^{*}$-algebra $A(\sum)=C(X)
\rtimes_{\alpha} \mathbb{Z}$. Thus the $*$-representations of the
$C^{*}$-crossed product $A(\sum)=C(X) \rtimes_{\alpha} \mathbb{Z}$
can be completely described and studied in terms of the covariant representation of the system 
$(C(X),\alpha , \mathbb{Z})$, that is
in terms of families of operators satisfying the covariance
commutation relations (\ref{covrelrep}) and the corresponding involution conditions. If $(\pi,u)$ is a
covariant representation of the system 
$(C(X),\alpha ,\mathbb{Z})$, then the corresponding $*$-representation of the crossed
product $C^{*}$-algebra $C(X) \rtimes_{\alpha} \mathbb{Z}$
transforms a generalized polynomial $\sum_{j=-n}^{n} f_j
\delta^{j}$ into the operator $\sum_{j=-n}^{n} \pi(f_j) u^{j}$.

\section{$\cs$-algebras of non-invertible dynamical systems}
The $\cs$-crossed product by $\Z$ is an important way to associate a
$\cs$-algebra to an invertible dynamical system. There are several
ways to generalize this construction to non-invertible dynamical
systems. One of these is due to Exel. It relies on transfer
operators.

We will here give a short description of Exel's construction:

\begin{definition}
  A $\cs$-dynamical system is a pair $(A,\alpha)$ of a unital
  $\cs$-algebra $A$ and an $*$-endomorphism $\alpha:A\to A$.
\end{definition}

\begin{definition}
  A transfer operator for the $\cs$-dynamical system $(A,\alpha)$ is a
  continuous linear map $\tra:A\to A$ such that
  \begin{enumerate}
  \item $\tra$ is positive in the sense that $\tra(A_+)\subseteq A_+$,
  \item $\tra(\alpha(a)b)=a\tra(b)$ for all $a,b\in A$.
  \end{enumerate}
\end{definition}

\begin{definition}
  Given a $\cs$-dynamical system $(A,\alpha)$ and a transfer operator
  $\tra$ of $(A,\alpha)$, we let $\toep(A,\alpha,\tra)$ be the
  universal unital $\cs$-algebra generated by a copy of $A$ and an
  element $S$ subject to the relations
  \begin{enumerate}
  \item $Sa=\alpha(a)S$,
  \item $S^*aS=\tra(a)$,
  \end{enumerate}
  for all $a\in A$.
\end{definition}

Using \cite{MR813640}, it is easy to see that relations are admissible
and thus that $\toep(A,\alpha,\tra)$ exists.

It is proved in \cite{MR2032486}*{Corollary 3.5} that the standard embedding of $A$ into $\toepl$ is injective. We will therefore
from now on view $A$ as a $\cs$-sub\-algebra of $\toepl$.

\begin{definition}
  By a redundancy we will mean a pair 
  $(a,k)\in A\times \overline{ASS^*A}$ such that $abS=kbS$ for all $b\in A$.
\end{definition}

\begin{definition}
  The crossed product $\cros$ is the quotient of $\toepl$ by the
  closed two-sided ideal generated by the set of differences $a-k$,
  for all redundancies $(a,k)$ such that $a\in \overline{A\alpha(A)A}$.
\end{definition}
We will denote the quotient map from $\toepl$ to $\cros$ by $\quo$.

If $(A,\alpha)$ is an invertible $\cs$-dynamical system, meaning
that $\alpha$ is an automorphism, then $\alpha\inv$ is a transfer
operator for $(A,\alpha)$.

Let us consider $\toep(A,\alpha,\alpha\inv)$. It follows from (2) that
$S^*S=I$, where $I$ denotes the unit of $A$. For all $b\in A$, 
\begin{equation*}
  SS^*bS=S\alpha\inv(b)=bS=IbS,
\end{equation*}
so $(I,SS^*)$ is a redundancy. Thus $\quo(S)$ is a unitary which
satisfies
\begin{equation*}
  \quo(S)\quo(a)\quo(S)^*=\quo(\alpha(a))
\end{equation*}
for all $a\in A$. In other words, $(\quo,\quo(S))$ is a covariant representation of
$(A,\Z,\alpha)$.

On the other hand, in $A\rtimes_\alpha\Z$, $\delta_1$ satisfies
\begin{enumerate}
\item $\delta_1a=\alpha(a)\delta_1$,
\item $\delta_1^*a\delta_1=\alpha\inv(a)$,
\end{enumerate}
for all $a\in A$, and if $ab\delta_1=kb\delta_1$ for all $b\in A$,
then
\begin{equation*}
  a=aI=aI\delta_1\delta_1^*=kI\delta_1\delta_1^*=kI=k,
\end{equation*}
so $a-k=0$ for all redundancies $(a,k)$, and thus $A\rtimes_\alpha\Z$
is isomorphic to $\cro{A}{\alpha}{\alpha\inv}$.

\section{Shift spaces} \label{notation}
For an introduction to shift spaces see \cite{MR97a:58050}.

Let $\al$ be a finite set endowed with the discrete topology. We will
call this set the alphabet. Let $\al^\No$ be the infinite
product spaces $\prod_{n=0}^\infty \al$
endowed with the product topology. The transformation
$\tsh$ on $\al^\No$ given by $(\tsh(x))_i=x_{i+1},\ i\in
\No$ is called the \emph{shift}. Let $\OSS$ be a shift invariant
closed subset of $\al^\No$ (by shift invariant we mean that
$\tsh(\OSS)\subseteq \OSS$, not necessarily
$\tsh(\OSS)=\OSS$). The topological dynamical system
$(\OSS,\tsh_{|\OSS})$ is called a \emph{shift space}. We will denote
$\tsh_{|\OSS}$ by $\tsh_{\OSS}$ or $\tsh$ for simplicity,
and on occasion the alphabet $\al$ by $\al_{\OSS}$.

We denote the $n$-fold composition of $\osh$ with itself by
$\osh^n$, and we denote the preimage of a set $X$ under $\osh^n$ by
$\osh^{-n}(X)$.


A finite sequence $u=(u_1,\ldots ,u_k)$ of elements $u_i\in
\al$ is called a finite \emph{word}. The \emph{length} of $u$ is $k$
and is denoted by $|u|$. We let for each $k\in \No$, $\al^k$ be the
set of all words with length $k$ and we let
$\B^k(\OSS)$ be the set of all words with length $k$ appearing in
some $x\in \OSS$. We set $\B_l(\OSS)=\bigcup_{k=0}^l
\B^k(\OSS)$ and $\B(\OSS)=\bigcup_{k=0}^\infty \B^k(\OSS)$ and likewise
$\al_l=\bigcup_{k=0}^l\al^k$ and $\alwords=\bigcup_{k=0}^\infty \al^k$,
where $\B^0(\OSS)=\al^0$ denote the set consisting of the empty word
$\emptyword$ which has length 0. $\B(\OSS)$ is called the
\emph{language} of $\OSS$. Note
that $\B(\OSS)\subseteq \alwords$ for every shift space.

For a shift space $\OSS$ and a word $u \in \B(\OSS)$ we denote by
$C_{\OSS}(u)$ the \emph{cylinder set}
\[C_{\OSS}(u)=\{x\in \OSS \mid (x_1,x_2,\dots
,x_{|u|})=u\}.\]
It is easy to see that
\[\{C_{\OSS}(u)\mid u \in \B(\OSS)\}\]
is a basis for the topology of $\OSS$, and that $C_{\OSS}(u)$ is
closed and compact for every $u\in \B(\OSS)$. We will allow
ourself to write $C(u)$ instead of $C_{\OSS}(u)$ when it is clear
which shift space space we are working with.

For a shift space $\OSS$ and words $u, v \in \B(\OSS)$ we
denote by $\cy{u}{v}$ the set
\[C(v)\cap\tsh^{-|v|}(\tsh^{|u|}(C(u)))=\{v x\in \OSS\mid u x\in \OSS\}.\]

What we have defined above is a \emph{one-sided} shift space. A
\emph{two-sided} shift space is defined in the same way, except that we
replace $\No$ with $\Z$: Let $\al^\Z$ be the infinite
product spaces $\prod_{n=-\infty}^\infty \al$
endowed with the product topology, and let $\tsh$ be the transformation
on $\al^\Z$ given by $(\tsh(x))_i=x_{i+1},\ i\in
\Z$. A shift invariant closed subset
$\Lambda$ of $\al^\Z$  (here, by shift
invariant we mean $\tsh(\Lambda)=\Lambda$) is called a \emph{two-sided
  shift space}. The set
\begin{equation*}
  \OSS_\Lambda=\{(x_i)_{i\in\No}\mid (x_i)_{i\in\Z}\in \Lambda\}
\end{equation*}
is a one-sided shift space, and it is called the one-sided shift space of
$\Lambda$.

If $\OSS$ and $\OSSY$ are two shift spaces and $\phi:\OSS\to\OSSY$ is a
homeomorphism such that
$\phi\circ{\osh}_{\OSS}={\osh}_{\OSSY}\circ\phi$, then we say that $\phi$
is a \emph{conjugacy} and that $\OSS$ and $\OSSY$ are \emph{conjugate}
or \emph{one-sided conjugate} if we want to emphasis that we are
dealing with one-sided shift spaces. Likewise we say that two
two-sided shift spaces $\Lambda$ and $\Gamma$ are \emph{two-sided
  conjugate} if there exists a homeomorphism $\phi:\Lambda\to\Gamma$
such that $\phi\circ{\osh}_{\Lambda}={\osh}_{\Gamma}\circ\phi$. It is
an easy exercise to prove that if $\OSS_\Lambda$ and $\OSS_\Gamma$ are
one-sided conjugate, then $\Lambda$ and $\Gamma$ are two-sided
conjugate.

The weaker notion of \emph{flow equivalence} among two-sided shift spaces is
also of importance here.
This notion is defined using the
\emph{suspension flow space} of $(\Lambda,\tsh)$ defined as
$S\Lambda=(\Lambda\times\R)/\!\sim$ where  the equivalence relation
$\sim$ is generated by requiring that
$(x,t+1)\sim(\tsh(x),t)$. Equipped with the quotient topology, we get
a compact space with a \emph{continuous flow} consisting of a family of
maps $(\phi_t)$ defined by $\phi_t([x,s])=[x,s+t]$. We say that two
two-sided shift spaces $\Lambda$ and $\Gamma$ are
\emph{flow equivalent} and write $\Lambda\floweq\Gamma$ if a
homeomorphism $F:S\Lambda\to S\Gamma$ exists with the property that
for every $x\in S\Lambda$ there is a monotonically increasing map
$f_x:\R\to\R$ such that
\[
F(\phi_t(x))=\phi'_{f_x(t)}(F(x)).
\]
In words, $F$ takes flow orbits to flow orbits in an
orientation-preserving way.
It is not hard to see that two-sided conjugacy implies flow equivalence.

\section{The $\cs$-algebra associated with a shift space}

Let $(\OSS,\osh)$ be a one-sided shift space. We want to define on  
$C(\OSS)$ a
transfer operator 
\begin{equation*}
  \tra(f)(x)=
  \begin{cases}
    \frac{1}{\#\osh\inv(\{x\})}\sum_{y\in\osh\inv\{x\}}f(y)&\text{if
    }x\in \osh(\OSS),\\
    0& \text{if }x\notin\osh(\OSS),
  \end{cases}
\end{equation*}
where the symbol $\# $ is used for the cardinality of a set. But such operator might take us out of the class of continuous functions on
$\OSS$. So we let $\DX$ be the smallest $\cs$-subalgebra of the
$\cs$-algebra of bounded functions on $\OSS$, containing $C(\OSS)$ and
closed under $\tra$ and $\alpha$, where $\alpha$ is the map $f\mapsto
f\circ\osh$.

\begin{lemma} \label{lemma:number}
  The function
  \begin{equation*}
    x\mapsto \osh^{-n}\{x\},\ x\in\OSS
  \end{equation*}
  belongs to $\DX$ for every $n\in\N$.
\end{lemma}

\begin{proof}
  Let $n\in\N$.
  The function
  \begin{equation*}
    f=1-\tra^n(1)+\sum_{u\in\al^n}(\tra^n(1_{\cyl{u}}))^2
  \end{equation*}
  belongs to $\DX$, and since
  \begin{equation*}
    f(x)=
    \begin{cases}
      \frac{1}{\#\osh^{-n}\{x\}}&\text{if }x\in\osh^n(\OSS),\\
      1&\text{if }x\notin\osh^n(\OSS),
    \end{cases}
  \end{equation*}
  for every $x\in\OSS$, $f$ is invertible and $f\inv\in\DX$, and so does
  $f\inv+\tra^n(1)-1$. Since
  \begin{equation*}
    (f\inv+\tra^n(1)-1)(x)=
    \begin{cases}
      \#\osh^{-n}\{x\}&\text{if }x\in\osh^n(\OSS),\\
      0&\text{if }x\notin\osh^n(\OSS),
    \end{cases}
    =\#\osh^{-n}\{x\}
  \end{equation*}
  for every $x\in \OSS$, we are done.
\end{proof}

\begin{lemma} \label{lemma:generated}
  $\DX$ is the $\cs$-algebra generated by $\{1_{\cy{u}{v}}\}_{u,v\in
    \alwords}$.
\end{lemma}

\begin{proof}
  Let $f(x)=\#\osh^{-|u|}\{x\}$. It then follows from
  Lemma \ref{lemma:number} that $f\in\DX$. Thus
  \begin{equation*}
    1_{\cy{u}{v}}=1_{\cyl{v}}\alpha^{|v|}(f\tra^{|u|}(1_{\cyl{u}}))
  \end{equation*}
  belongs to $\DX$ for every $u,v\in \alwords$, and $\cs(1_{\cy{u}{v}}\mid
  u,v\in\alwords)\subseteq \DX$.

  In the other direction we have that since $\{\cyl{v}\}_{v\in\alwords}$ is
  a basis of the topology of $\OSS$ consisting of clopen sets,
  $\{1_{\cyl{v}}\}_{v\in\alwords}$ generates $C(\OSS)$ and since
  $1_{\cyl{v}}=1_{\cy{\emptyword}{v}}$, $C(\OSS)$ is
  contained in $\cs(1_{\cy{u}{v}}\mid u,v\in\alwords)$. Since
  \begin{equation*}
    \alpha(1_{\cy{u}{v}})=\sum_{a\in\al}1_{\cy{u}{av}},
  \end{equation*}
  and
  \begin{equation*}
      \tra(1_{\cy{u}{v}})= \left(\sum_{a\in\al}(\tra(1_{\cyl{a}}))^2\right)
      1_{\cy{v_1}{\emptyword}} 1_{\cy{u}{v_2v_3\dotsm v_n}},
    \end{equation*}
    if $v\ne \emptyword$, and
    \begin{equation*}
       \tra(1_{\cy{u}{\emptyword}})=
       \left(\sum_{a\in\al}(\tra(1_{\cyl{a}}))^2\right)
       \left(\sum_{a\in\al} 1_{\cy{au}{\emptyword}}\right),
    \end{equation*}
    the $\cs$-algebra generated by $1_{\cy{u}{v}}$, $u,v\in \alwords$ is
    closed under $\tra$ and $\alpha$ and thus contain $\DX$.
\end{proof}

\begin{theorem} \label{universal}
  The $\cs$-algebra $\crod$ is the universal
  $\cs$-algebra generated by partial isometries $\{S_u\}_{u\in\alwords}$
  satisfying:
  \begin{enumerate}
  \item $S_uS_v=S_{uv}$ for all $u,v\in\alwords$,
  \item the map
    \begin{equation*}
      1_{\cy{u}{v}}\mapsto S_vS_u^*S_uS_v^*,\ u,v\in\alwords
    \end{equation*}
    extends to a $*$-homomorphism from $\DX$ to $\cs(S_u\mid u\in\alwords)$.
  \end{enumerate}
\end{theorem}

\begin{proof}
  We will first show that $\crod$ is generated by
  partial isometries $\{S_u\}_{u\in\alwords}$ satisfying $(1)$ and $(2)$,
  and then that if $A$ is a $\cs$-algebra generated by partial
  isometries $\{s_u\}_{u\in\alwords}$ satisfying $(1)$ and $(2)$, then
  there is a $*$-homomorphism from $\crod$ to $A$
  sending $S_u$ to $s_u$ for all $u\in\alwords$.

  Working within $\toep(\DX,\alpha,\tra)$ let for each $a\in\al$,
  \begin{equation*}
    T_a=1_{\cyl{a}}\left(\alpha(f)\right)^{1/2}S,
  \end{equation*}
  where $f$ is the function $x\mapsto \#\osh\inv\{x\}$, which belongs
  to $\DX$ by Lemma \ref{lemma:number}, and we let for each
  $u=u_1u_2,\dotsm u_n\in\alwords$,
  \begin{equation*}
    S_u=\quo(T_{u_1})\quo(T_{u_2})\dotsm \quo(T_{u_n}).
  \end{equation*}
  Then clearly $\{S_u\}_{u\in\alwords}$ satisfy $(1)$.

  Let $a\in\al$, $g\in\DX$ and $x\in\OSS$. Then
  \begin{equation*}
    \begin{split}
      \left(T_a^*gT_a\right)(x)&= \left(S^* \alpha(f)1_{\cyl{a}}g S\right)(x)\\
      &= \left(\tra\left(\alpha(f)1_{\cyl{a}}g\right)\right)(x)\\
      &= \left(f \tra\left(1_{\cyl{a}}g\right)\right)(x)\\
      &=
      \begin{cases}
        g(ax)& \text{if }ax\in \OSS,\\
        0 & \text{if }ax\notin \OSS.
      \end{cases}
    \end{split}
  \end{equation*}

  Now let $a\in\al$ and $g,h\in\DX$. Then
  \begin{equation*}
    \begin{split}
      T_agT_a^*hS &=
      1_{\cyl{a}}\left(\alpha(f)\right)^{1/2}SgS^*\left(\alpha(f)\right)^{1/2}1_{\cyl{a}}hS\\
      &=1_{\cyl{a}}\left(\alpha(f)\right)^{1/2}Sg\tra\left(\left(\alpha(f)\right)^{1/2}1_{\cyl{a}}h\right)\\
       &=1_{\cyl{a}}\alpha\left(f^{1/2}g\tra\left(\left(\alpha(f)\right)^{1/2}1_{\cyl{a}}h\right)\right)S,
    \end{split}
  \end{equation*}
  and for every $x\in\OSS$ is
  \begin{equation*}
    \left(1_{\cyl{a}}\alpha\left(f^{1/2}g\tra\left(\left(\alpha(f)\right)^{1/2}1_{\cyl{a}}h\right)\right)\right)(x)=
    \begin{cases}
      g(\osh(x))h(x)&\text{if }x\in\cyl{a},\\
      0&\text{if }x\notin\cyl{a}.
    \end{cases}
  \end{equation*}
  Since $SgS^*=\alpha(g)SS^*$,
  \begin{equation*}
     T_agT_a^*=
     1_{\cyl{a}}\left(\alpha(f)\right)^{1/2}SgS^*\left(\alpha(f)\right)^{1/2}1_{\cyl{a}} \in \overline{ASS^*A},
  \end{equation*}
  so $(\alpha(g)1_{\cyl{a}}, T_agT_a^*)$ is a redundancy, and since
  $\alpha(g)1_{\cyl{a}}\in \overline{A\alpha(A)A}$, it follows
  that $S_a\quo(g)S_a^*$ and $\quo(\alpha(g)1_{\cyl{a}})$ are equal in
  $\crod$.

  Thus $S_a^*\quo(g)S_a=\lamb{a}(g)$ and
  $S_a\quo(g)S_a^*=\quo(\alpha(g)1_{\cyl{a}})$ for every $a\in\al$ and
  $g\in \DX$, where $\lamb{a}(g)$ is the map given by
  \begin{equation*}
    \lamb{a}(g)(x)=
      \begin{cases}
          g(ax)& \text{if }ax\in \OSS,\\
          0 & \text{if }ax\notin \OSS,
    \end{cases}
  \end{equation*}
  for $x\in\OSS$, which shows that
  \begin{equation*}
    \quo(1_{\cy{u}{v}})=S_vS_u^*S_uS_v^*
  \end{equation*}
  for every $u,v\in\alwords$. Hence $\{S_u\}_{u\in\alwords}$ satisfy $(2)$. To
  see that $\crod$ is generated by $\{S_u\}_{u\in\alwords}$, we first notice
  that $\toep(\DX,\alpha,\tra)$ is generated by $\DX$ and $S$, and
  that $\DX$, by Lemma \ref{lemma:generated}, is generated by
  $\{1_{\cy{u}{v}}\}_{u,v\in \alwords}$, and then that the function
  $\alpha(f)$, where $f$ as before is the function
  \begin{equation*}
    x\mapsto \osh^{-n}\{x\},\ x\in\OSS,
  \end{equation*}
  is invertible and that
  $S=\sum_{a\in\al}\alpha(f)^{-1/2}T_a$. Thus it follows that $\crod$
  is generated by $\{S_u\}_{u\in\alwords}$.

  Assume now that $A$ is a $\cs$-algebra generated by partial
  isometries $\{s_u\}_{u\in\alwords}$ which satisfy $(1)$ and $(2)$. We
  then let
  $s=\sum_{a\in\al}\phi(\alpha(f)^{-1/2})s_a$, where $\phi$
  is the $*$-homomorphism from $\DX$ to $\cs(s_u\mid u\in\alwords)$,
  which extends the map
  \begin{equation*}
    1_{\cy{u}{v}}\mapsto s_vs_u^*s_us_v^*,\ u,v\in\alwords.
  \end{equation*}
  Observe first that if $a,b\in\al$ and $a\ne b$, then $s_a^*s_b=0$,
  because
  \begin{equation*}
    s_a^*s_b=s_a^*s_as_a^*s_bs_b^*s_b=s_a^*\phi(1_{\cyl{a}}1_{\cyl{b}})s_b.
  \end{equation*}

  Let then $a\in\al$ and $u,v\in\alwords$. We then have that if $v\ne
  \emptyword$, then
  \begin{equation*}
    \begin{split}
      s_a^* \phi(1_{\cy{u}{v}}) s_a &= s_a^* s_v s_u^* s_u s_v^* s_a\\
      &= s_a^* s_{v_1} s_{v_2v_3\dotsm v_{|v|}} s_u^* s_u
      s_{v_2v_3\dotsm v_{|v|}}^* s_{v_1}^* s_a\\
      &=
      \begin{cases}
        \phi(1_{\cy{v_1}{\emptyword}}1_{\cy{u}{v_2v_3\dotsm
            v_{|v|}}})& \text{if }a=v_1,\\
        0& \text{if }a\ne v_1
      \end{cases}\\
      &= \phi(\lamb{a}(1_{\cy{u}{v}})),
    \end{split}
  \end{equation*}
  and if $v=\emptyword$, then
  \begin{equation*}
    \begin{split}
       s_a^* \phi(1_{\cy{u}{v}}) s_a &= s_a^* s_u^* s_u s_a\\
       &= s_{ua}^* s_{ua}\\
       &= \phi(1_{\cy{ua}{\emptyword}})\\
       &= \phi(\lamb{a}(1_{\cy{u}{v}})).
    \end{split}
  \end{equation*}
  Since $\DX$ is generated by $1_{\cy{u}{v}}$, $u,v\in \alwords$, this
  shows that $s_a^*\phi(g)s_a=\phi(\lamb{a}(g))$ for each $a\in\al$
  and every $g\in\DX$.

    Let $u,v\in\alwords$. Then
  \begin{equation*}
    \begin{split}
      \phi(\alpha(1_{\cy{u}{v}})) s &= \sum_{a\in\al} \phi(1_{\cy{u}{av}})s\\
        &= \sum_{a\in\al} \phi(1_{\cy{u}{av}}) \sum_{b\in\al}
        \phi(\alpha(f)^{-1/2}) s_b\\
        &= \sum_{a\in\al} \sum_{b\in\al} \phi(\alpha(f)^{-1/2}
        1_{\cy{u}{av}}) s_bs_b^*s_b\\
        &= \sum_{a\in\al} \sum_{b\in\al} \phi(\alpha(f)^{-1/2}
        1_{\cy{u}{av}} 1_{\cyl{b}}) s_b\\
        &= \sum_{a\in\al} \phi(\alpha(f)^{-1/2}1_{\cy{u}{av}}) s_a\\
        &= \sum_{a\in\al} \phi(\alpha(f)^{-1/2})s_{av}s_u^*s_us_{av}^* s_a\\
        &= \sum_{a\in\al} \phi(\alpha(f)^{-1/2}) s_as_vs_u^*s_us_v^*s_a^*s_a\\
        &= \sum_{a\in\al} \phi(\alpha(f)^{-1/2}) s_a \phi(1_{\cy{u}{v}}
        1_{\cy{a}{\emptyword}})\\
        &= \sum_{a\in\al} \phi(\alpha(f)^{-1/2}) s_a
        \phi(1_{\cy{a}{\emptyword}} 1_{\cy{u}{v}})\\
        &= \sum_{a\in\al} \phi(\alpha(f)^{-1/2}) s_a s_a^* s_a
        \phi(1_{\cy{u}{v}})\\
        &= \sum_{a\in\al} \phi(\alpha(f)^{-1/2}) s_a \phi(1_{\cy{u}{v}})\\
        &= s \phi(1_{\cy{u}{v}}),
    \end{split}
  \end{equation*}
  and
  \begin{equation*}
    \begin{split}
      s^*\phi(1_{\cy{u}{v}})s &= \sum_{a\in\al} \sum_{b\in\al} s_a^*
      \phi(\alpha(f)^{-1/2} 1_{\cy{u}{v}} \alpha(f)^{-1/2}) s_b\\
      &= \sum_{a\in\al} \sum_{b\in\al} s_a^* \phi(\alpha(f)^{-1}
      1_{\cyl{a}} 1_{\cyl{b}} 1_{\cy{u}{v}}) s_b\\
      &= \phi(\alpha(f)^{-1}) \sum_{a\in\al} s_a^* \phi(1_{\cy{u}{v}})
      s_a^*\\
      &= \phi(\alpha(f)^{-1}) \sum_{a\in\al} s_a^* s_vs_u^*s_us_v^*
      s_a^*\\
      &= \phi(\alpha(f)^{-1}) \sum_{a\in\al} s_a^*
      s_{v_1}s_{v_2v_3\dotsm v_{|v|}} s_u^*s_u s_{v_2v_3\dotsm
        v_{|v|}}^* s_{v_1}^* s_a^*\\
      &= \phi(\alpha(f)^{-1}) s_{v_1}^* s_{v_1} s_{v_2v_3\dotsm
        v_{|v|}} s_u^*s_u s_{v_2v_3\dotsm
        v_{|v|}}^* s_{v_1}^* s_{v_1}^*\\
      &= \phi(\alpha(f)^{-1} 1_{\cy{v_1}{\emptyword}} 1_{\cy{v_2v_3\dotsm
        v_{|v|}}{u}})\\
      &= \phi(\tra(1_{\cy{v}{u}}))
    \end{split}
  \end{equation*}
  and since $\DX$ is generated by $1_{\cy{u}{v}}$, $u,v\in \alwords$, this
  shows that $\phi(\alpha(g))s=s\phi(g)$ for every $g\in\DX$.

  Thus it follows from the universal property of $\toepd$, that there
  exists a $*$-homomorphism $\psi$ from $\toepd$ to $A$ which maps $g$
  to $\phi(g)$ for $g\in \DX$ and $S$ to $s$. We will now show that
  $\psi$ vanishes on the closed two-sided ideal generated by the set of
  differences $g-k$, for all redundancies $(g,k)$ such that $g\in
  \overline{\DX\alpha(\DX)\DX}$, and thus that it factors through the
  quotient and yields a $*$-homomorphism $\widetilde{\psi}:\crod\to A$
  such that $\widetilde{\psi}(\quo(g))=\phi(g)$ and
  $\widetilde{\psi}(\quo(S))=s$, and hence
  \begin{equation*}
    \begin{split}
      \widetilde{\psi}(S_a) &= \widetilde{\psi}(\quo(T_a))\\
      &= \widetilde{\psi}(\quo(1_{\cyl{a}}(\alpha(f)^{1/2})S))\\
      &= \phi(1_{\cyl{a}}(\alpha(f))^{1/2})s\\
      &= \phi(1_{\cyl{a}}(\alpha(f))^{1/2}) \sum_{b\in\al}
      \phi((\alpha(f))^{-1/2})s_b\\
      &= \sum_{b\in\al}\phi(1_{\cyl{a}}1_{\cyl{b}})s_b\\
      &= \phi(1_{\cyl{a}})s_a\\
      &= s_as_a^*s_a\\
      &= s_a
    \end{split}
  \end{equation*}
  for all $a\in\al$, and thus $\widetilde{\psi}(S_u)=s_u$ for every
  $u\in\alwords$.

  Assume that $g\in\overline{\DX\alpha(\DX)\DX}$, that $k\in
  \overline{\DX SS^*\DX}$ and $ghS=khS$ for every $h\in \DX$. Then
  \begin{equation*}
    \begin{split}
      \psi(g) &= \psi(g\sum_{a\in\al}1_{\cyl{a}})= \psi(g)\sum_{a\in\al}s_as_a^*\\
      &= \psi(g) \sum_{a\in\al} \phi(1_{\cyl{a}}(\alpha(f))^{1/2})ss_a^*\\
      &= \sum_{a\in\al} \psi(g1_{\cyl{a}}(\alpha(f))^{1/2}S)s_a^*\\
      &= \sum_{a\in\al} \psi(k1_{\cyl{a}}(\alpha(f))^{1/2}S)s_a^*\\
      &= \psi(k) \sum_{a\in\al}
      \phi(1_{\cyl{a}}(\alpha(f))^{1/2})ss_a^*\\
      &= \psi(k)\sum_{a\in\al}s_as_a^* = \psi(k\sum_{a\in\al}1_{\cyl{a}})\\
      &= \psi(k),
    \end{split}
  \end{equation*}
  so $\psi$ vanishes on the closed two-sided ideal generated by the set of
  differences $g-k$, for all redundancies $(g,k)$ such that $g\in
  \overline{\DX\alpha(\DX)\DX}$.
\end{proof}

\section{A representation of $\crod$} \label{representation}
Let $\OSS$ be a shift space, let $\HH_{\OSS}$ be the Hilbert space
$l_2(\OSS)$ and let $\{e_x\}_{x\in\OSS}$ be an orthonormal basis for
$\HH_{\OSS}$. Let for every $u\in\alwords$, $s_u$ be the operator on
$\HH_{\OSS}$ defined by
\begin{equation*}
  s_u(e_x)=
  \begin{cases}
    e_{ux}& \text{if }ux\in\OSS,\\
    0& \text{if }ux\notin\OSS.
  \end{cases}
\end{equation*}
We leave it to the reader to check that the operators $\{s_u\}_{u\in\alwords}$
satisfy condition (1) and (2) of Theorem \ref{universal}. Thus there
exists a $*$-homomorphism $\phi:\crod\to \cs(s_u\mid u\in\alwords)$
such that $\phi(S_u)=s_u$ for every $u\in\alwords$. In other words,
$S_u\mapsto s_u$ is a representation of $\crod$ on the Hilbert space
$\HH_{\OSS}$.

This representation is in general not faithful. If for example $\OSS$
only consist of one word, then $\crod$ is isomorphic to $C(\T)$,
whereas $\cs(s_u\mid u\in\alwords)$ is isomorphic to $\C$. We will in
section \ref{sec:uniq} see that if the shift space $\OSS$ satisfies a certain
condition $(I)$, then the representation $\phi$ is injective. We will
in section \ref{sec:uniq} construct a representation of $\crod$ which
is faithful for every shift space $\OSS$.

Although the $*$-homomorphism $\phi:\crod\to \cs(s_u\mid
u\in\alwords)$ is not in general injective the restriction of $\phi$
to $\DX$ is, and so it follows from the universal property of $\crod$,
that also the restriction of $\rho:\toepd\to\crod$ to $\DX$ is
injective. Thus we will allow ourselves to view $\DX$ as a sub-algebra
of $\crod$. We then have 
\begin{equation*}
  1_{\cy{u}{v}}= S_vS_u^*S_uS_v^*
\end{equation*}
for all $u,v\in\alwords$.

\section{$\crod$'s relationship with other $\cs$-algebras associated
  to shift spaces}
As mentioned in the introduction, other $\cs$-algebras have been
associated to shift spaces.
We will in this section look at the relation between these
$\cs$-algebras and $\crod$.

As far as the authors know, three different construction of
$\cs$-algebras associated to shift spaces appears in the
literature. These are:
\begin{itemize}
\item The $\cs$-algebra $\Oo_\Lambda$ defined in \cite{MR98h:46077},
\item the $\cs$-algebra $\Oo_\Lambda$ defined in \cite{MR2091486},
\item the $\cs$-algebra $\Oo_{\OSS}$ defined in \cite{tmc}.
\end{itemize}
These are all $\cs$-algebras generated by partial isometries
$\{S_a\}_{a\in\al}$, where $\al$ is the alphabet of the shift space in
question. The two first $\cs$-algebras are defined for every two-sided
shift space $\Lambda$, whereas the last one is defined for every
one-sided shift space $\OSS$.

We will in this section see, that there for every one-sided shift
space $\OSS$ exists a $*$-isomorphism between $\crod$ and the
$\cs$-algebra $\Oo_{\OSS}$ defined in \cite{tmc} which maps $S_a$ to
$S_a$ for every $a\in\al$, and that there for every two-sided shift
space $\Lambda$ exist a surjective $*$-homomorphism from the
$\cs$-algebra $\Oo_\Lambda$ defined in \cite{MR98h:46077} to $\crol$
which maps $S_a$ to $S_a$ for every $a\in\al$, and a surjective
$*$-homomorphism from $\crol$ to the
$\cs$-algebra $\Oo_\Lambda$ defined in \cite{MR2091486}
which maps $S_a$ to $S_a$ for every $a\in\al$. The first of these
surjective $*$-homomorphisms is injective if $\Lambda$ satisfy the
condition (*) defined in \cite{MR2091486}, and the second surjective
$*$-homomorphism is injective if $\Lambda$ satisfy the condition $(I)$ in Section \ref{sec:uniq}.

\begin{remark} \label{corresp}
In \cite{tmc}, a $\cs$-algebra $\Oo_{\OSS}$ has been constructed by
using $\cs$-correspondences and Cuntz-Pimsner algebras for every shift
space $\OSS$. It follows from Theorem \ref{universal} and
\cite{tmc}*{Remark 7.4} that $\Oo_{\OSS}$ is isomorphic to $\crod$ for
every one-sided shift space $\OSS$. Thus it follows from
\cite{tmc}*{Remark 7.4} that for every two-sided shift space
$\Lambda$, the algebra $\crol$ satisfy all of the results the algebra
$\Oo_\Lambda$ is claimed to satisfy in
\cites{MR1911208,MR2002i:37015,MR2002h:19004,MR2001g:46147,MR2001e:46115,MR2000d:46082,MR2000f:46084,MR2000e:46087,MR98h:46077}.
\end{remark}

\begin{remark}
  In \cite{MR2091486} a $\cs$-algebra $\Oo_\Lambda$ has been defined
  for every two-sided shift space by defining operators on the Hilbert
  space $l_2(\OSS_\Lambda)$. These operators are identical to the
  operators $s_u$ defined in section \ref{representation} for $\OSS$
  equal to the one-sided shift space $\OSS_\Lambda$ associated to
  $\Lambda$. Thus we have for every two-sided shift space $\Lambda$ a
  surjective $*$-homomorphism from $\crol$ to $\Oo_\Lambda$ which is
  injective if $\Lambda$ satisfies condition $(I)$, and we
  also know that there are examples of two-sided shift spaces (for
  instance the shift only consisting of one point) for which the
  $*$-homomorphism is not injective.
\end{remark}

As we have mentioned before, our $\cs$-algebra $\crol$ satisfies all of the results that the algebra
$\Oo_\Lambda$ is claimed to satisfy
\cites{MR1911208,MR2002i:37015,MR2002h:19004,MR2001g:46147,MR2001e:46115,MR2000d:46082,MR2000f:46084,MR2000e:46087,MR98h:46077},
whereas the $\cs$-algebra $\Oo_\Lambda$ originally defined in \cite{MR2001g:46147}, does
not. The latter $\cs$-algebra have been properly characterized in
\cite{MR2091486} (where it is called $\Oo_\Lambda^*$). We will now use
this characterization to show that
there for every two-sided shift space $\Lambda$ exists a surjective
$*$-homomorphism from $\Oo_\Lambda$ to $\crol$.

Let for every $l\in\No$, $\A_l^*$ be the $\cs$-subalgebra of $\Oo_\Lambda$
generated by $\{S_u^*S_u\}_{u\in \al_l}$, and let $\A_\Lambda^*$ be
the $\cs$-subalgebra of $\Oo_\Lambda$ generated by $\{S_u^*S_u\}_{u\in
  \al}$. Notice that
\begin{equation*}
  \A_\Lambda^*=\overline{\bigcup_{l\in\No}\A_l^*}.
\end{equation*}
The key to characterizing $\Oo_\Lambda$ is to describe $\A_l^*$ and
$\A_\Lambda^*$. This is done in this way:

Let for every $l\in\No$ and every $u\in\B(\Lambda)$,
\begin{equation*}
  \Past_l(u)=\{v\in\al_l\mid vu\in\B(\Lambda)\}.
\end{equation*}
We then define an equivalence relation $\sim_l$ on $\B(\Lambda)$ called
\emph{$l$-past equivalence} in this way:
\begin{equation*}
  u\sim_l v\iff \Past_l(u)=\Past_l(v).
\end{equation*}
We denote the $l$-past equivalence class containing $u$ by $[u]_l$,
and we let
\begin{equation*}
  \B_l^*(\Lambda)=\{u\in\al_l^*\mid \text{the cardinality of
  }[u]_l\text{ is infinite}\},
\end{equation*}
and $\Omega_l^*=\B_l^*/\sim_l$. Since $\alwords_l$ is finite, so is
$\Omega_l^*$. We equip $\Omega_l^*$ with the discrete topology (so
$C(\Omega_l^*)\cong \C^{m^*(l)}$, where $m^*(l)$ is the number of elements of
$l$-past equivalence classes).

We then have:
\begin{lemma}[cf. \cite{MR2091486}*{Lemma 2.9}]
  The map
  \begin{equation*}
    1_{\{[u]_l\}}\mapsto 1_{[u]_l},\ u\in\B_l^*(\Lambda)
  \end{equation*}
  extends to a $*$-isomorphism between $C(\Omega_l^*)$ and $\A_l^*$.
\end{lemma}

We will now make the corresponding characterization of $\crod$:

Let $\OSS$ be a one-sided shift space. We let for every $l\in\No$,
$\A_l$ be the $\cs$-subalgebra of $\DX$ generated by
$\{1_{\cy{v}{\emptyword}}\}_{v\in\alwords_l}$, and we let $\AX$ be the
$\cs$-subalgebra of $\DX$ generated by
$\{1_{\cy{v}{\emptyword}}\}_{v\in\alwords}$. Notice that
\begin{equation*}
  \AX=\overline{\bigcup_{l\in\No}\A_l}.
\end{equation*}
Following Matsumoto (cf. \cite{MR2000d:46082}), we let for every
$l\in\N$ and every $x\in \OSS$,
\begin{equation*}
  \Past_l(x)=\{u\in\alwords_l\mid ux\in\OSS\}.
\end{equation*}
We then define an equivalence relation $\sim_l$ on $\OSS$ called
\emph{$l$-past equivalence} in this way:
\begin{equation*}
  x\sim_l y\iff \Past_l(x)=\Past_l(y).
\end{equation*}
We let $\Omega_l=\OSS/\sim_l$, and denote the $l$-past equivalence
class containing $x$ by $[x]_l$. Since $\alwords_l$ is finite, so is
$\Omega_l$. We equip $\Omega_l$ with the discrete topology (so
$C(\Omega_l)\cong \C^{m(l)}$, where $m(l)$ is the number of elements of
$l$-past equivalence classes). Since
\begin{equation*}
  [x]_l=\left(\bigcap_{u\in\Past_l(x)}\cy{u}{\emptyword}\right)\cap
  \left(\bigcap_{v\in\alwords_l\setminus\Past_l(x)}\OSS\setminus\cy{v}{\emptyword}\right),
\end{equation*}
the function $1_{[x]_l}$ belongs to $\A_l$, and
$\{1_{[x]_l}\}_{x\in\OSS}$ generates $\A_l$. Thus
\begin{equation*}
  1_{\{[x]_l\}}\mapsto 1_{[x]_l}
\end{equation*}
is a $*$-isomorphism between $C(\Omega_l)$ and $\A_l$, which extends
to an isomorphism between $C(\Omega_{\OSS})$ and $\AX$.

Consider the condition:

$(*):$ There exists for each $l \in \No$ and each infinite sequence of admissible words
$\{u_i\}_{i\in\N}$ satisfying $\Past_l(u_i) = \Past_l(u_j)$ for all
$i,j\in\N$, an $x\in\OSS_\Lambda$ such that
\begin{equation*}
\Past_l(x) = \Past_l(u_i)
\end{equation*}
for all $i\in\N$.

It follows from \cite{MR2091486}*{Corollary 3.3} that there is a
surjective $*$-homomorphism from $\A_\Lambda^*$ to $\A_{\OSS_\Lambda}$, and
that this $*$-homomorphism is injective if and only if $\Lambda$
satisfies the condition (*). As a consequence of this, we get that 
for every two-sided shift space $\Lambda$ exists a surjective
$*$-homomorphism from $\Oo_\Lambda$ to $\crol$, and that this
$*$-homomorphism is injective if $\Lambda$ satisfies the condition
(*).

There is in \cite{MR2091486} an example of a sofic shift space
$\Lambda$ for which $\Oo_\Lambda$ and $\crol$ are not isomorphic.

\section{Generalization of the Cuntz-Krieger algebras}
We are now able to show that $\crod$ in fact is a generalization
of the Cuntz-Krieger algebras. Actual we will prove that $\crod$
is a generalization of the universal Cuntz-Krieger algebra $\A\Oo_A$
that An Huef and Raeburn have constructed in \cite{MR98k:46098}.

\begin{theorem}
  Let $A$ be a $n\times n$-matrix with entries in $\{0,1\}$ and no zero
  rows, and let $\OSS_A$ be the one-sided shift spaces
  \begin{equation*}
    \left\{(x_i)_{i\in\No}\in \{1,2,\dotsc ,n\}^\No\mid \forall i\in
      \No: A(x_i,x_{i+1})=1\right\}.
  \end{equation*}
  Then $\croa$ is generated by partial isometries $\{S_i\}_{i\in
    \{1,2,\dotsc n\}}$ that satisfy
  \begin{equation*}
    \sum_{j=1}^nS_jS_j^*=I,
  \end{equation*}
  and
  \begin{equation*}
    S_i^*S_i=\sum_{j=1}^nA(i,j)S_jS_j^*
  \end{equation*}
  for every $i\in \{1,2,\dotsc ,n\}$.

  If $X$ is a unital $C^*$-algebra such that there exists a set of partial
  isometries $\{T_i\}_{i\in \{1,2,\dots ,n\}}$ in $X$ that satisfy
  \[\sum_{j=1}^nT_jT_j^*=I,\]
  and
  \[T_i^*T_i=\sum_{j=1}^nA(i,j)T_jT_j^*\]
  for every $i\in \{1,2,\dotsc ,n\}$; then there exists a
  $*$-homomorphism form $\croa$ to $X$ sending $S_i$ to
  $T_i$ for every $i\in \{1,2,\dotsc ,n\}$.
\end{theorem}

\begin{proof}
  Since $\OSS_A$ is the disjoint union of $\cyl{j},\ j\in\{1,2,\dotsc
  ,\}$,
  \begin{equation*}
    \sum_{j=1}^nS_jS_j^*=I,
  \end{equation*}
  and since for every $i\in\{1,2,\dotsc ,\}$, $\cy{i}{\emptyword}$ is
  the disjoint union of those $\cyl{j}$'s, where $A(i,j)=1$,
  \begin{equation*}
    S_i^*S_i=1_{\cy{i}{\emptyword}}=\sum_{j=1}^nA(i,j)1_{\cyl{j}}=
    \sum_{j=1}^nA(i,j)S_jS_j^*.
  \end{equation*}
The $\cs$-algebra $\croa$ is generated by partial isometries
  $\{S_u\}_{u\in\{1,2,\dotsc n\}^*}$, but since these partial isometries
  satisfy $S_uS_v=S_{uv}$ for all $u,v\in\{1,2,\dotsc n\}^*$,
  $\{S_i\}_{i\in\{1,2,\dotsc n\}}$ generates the whole $\croa$.

  Let $X$ be a unital $C^*$-algebra such that there exist partial
  isometries $T_i,\ i\in \{1,2,\dots ,n\}$ in $X$ that satisfy
  \[\sum_{j=1}^nT_jT_j^*=I,\]
  and
  \[T_i^*T_i=\sum_{j=1}^nA(i,j)T_jT_j^*\]
  for every $i\in \{1,2,\dots ,n\}$.  We let $T_\emptyword=I$ and we
  let for every $u=u_1u_2\dotsm
  u_n\in\{1,2,\dotsc ,n\}^*\setminus\{\emptyword\}$, $T_u$ be
  $T_u=T_{u_1}T_{u_2}\dotsm T_{u_n}$, and we will then show that
  \begin{enumerate}
  \item $T_uT_v=T_{uv}$ for all $u,v\in\{1,2,\dotsc,n\}^*$,
  \item the map
    \begin{equation*}
      1_{\cy{u}{v}}\mapsto T_vT_u^*T_uT_v^*,\ u,v\in\alwords
    \end{equation*}
    extends to a $*$-homomorphism from $\DX$ to $X$,
  \end{enumerate}
  and thus that there exists a $*$-homomorphism form $\croa$ to $X$ sending $S_u$ to
  $T_u$ for every $u\in \{1,2,\dotsc ,n\}^*$, and especially $S_i$ to
  $T_i$ for every $i\in \{1,2,\dotsc ,n\}$.

  It is clear from the way we defined $T_u$ that condition (1) is
  satisfied.
Let $m\in\N$, and denote by $\D_m$ the $C^*$-subalgebra of
  $\D_{\OSS_A}$ generated by $\{1_{\cyl{u}}\}_{u\in \{1,2,\dotsc
    ,n\}^m}$. If $u,v\in\{1,2,\dotsc ,n\}^m$ and $u\ne v$, then
  \begin{equation*}
    T_uT_u^*+T_vT_v^*
    \le \sum_{w\in\{1,2,\dotsc,n\}^m}T_wT_w^*=I,
  \end{equation*}
  and so
  \begin{equation*}
    T_u^* T_u+ T_u^* T_vT_v^* T_u=
    T_u^*(T_uT_u^*+T_vT_v^*)T_u \le T_u^*IT_u=T_u^*T_u,
  \end{equation*}
  which implies that $T_uT_u^*T_vT_v^*=T_uT_u^* T_vT_v^* T_uT_u^*=0$.

  Thus $\left\{T_uT_u^*\right\}_{u\in \{1,2,\dotsc ,n\}^m}$ are mutual
  orthogonal projections, and since $\left\{1_{\cyl{u}}\right\}_{u\in
    \{1,2,\dotsc ,n\}^m}$ also are mutual orthogonal
  projections and
  \begin{equation*}
    \begin{split}
      1_{\cyl{u}}=0 \ \Rightarrow\ & \cyl{u}=\emptyset\\
      \ \Rightarrow\ & u\notin \B(\OSS_A)\\
      \ \Rightarrow\ & \exists i\in\{1,2,\dotsc,m-1\}:A(u_i,u_{i+1})=0\\
      \ \Rightarrow\ & T_{u_i}T_{u_{i+1}} =
      T_{u_i}T_{u_i}^*T_{u_i}T_{u_{i+1}}T_{u_{i+1}}^*T_{u_{i+1}} =\\
      & T_{u_i}\sum_{k=1}^nA(U_i,k)T_kT_k^*T_{u_{i+1}}T_{u_{i+1}}^*T_{u_{i+1}} =
      0\\
      \ \Rightarrow\ & T_uT_u^*=0,
    \end{split}
  \end{equation*}
  there is a unital $*$-homomorphism $\psi_m$ from $\D_m$ to $X$
  obeying $\psi_m\left(1_{\cyl{u}}\right)=T_uT_u^*$ for every $u\in
  \{1,2,\dotsc ,n\}^m$.

  Since $\cyl{u}$ is the disjoint union of
  $\{\cyl{ui}\}_{i\in\{1,2,\dotsc,n\}}$,
  \begin{equation*}
    1_{\cyl{u}}=\sum_{i=1}^n1_{\cyl{ui}}\in D_{m+1}
  \end{equation*}
  for every $u\in\{1,2,\dotsc,n\}^m$, so $\D_m\subseteq \D_{m+1}$. Let
  us denote the inclusion of $\D_m$ into $\D_{m+1}$ by
  $\iota_m$. Since
  \begin{equation*}
    \begin{split}
      \psi_{m+1}\left(1_{\cyl{u}}\right)&=
      \psi_{m+1}\left(\sum_{i=1}^n1_{\cyl{ui}}\right)\\
      &= \sum_{i=1}^nT_{ui}T_{ui}^*\\
      &= T_u\left(\sum_{i=1}^nT_iT_i^*\right)T_u^*\\
      &= T_uT_u^*= \psi_m\left(1_{\cyl{u}}\right),
    \end{split}
  \end{equation*}
  $\psi_{m+1}\circ\iota_m=\psi_m$. Thus the $\psi_m$'s extends to a
  $*$-homomorphism from $\overline{\bigcup_{m\in\N}\D_m}$ to $X$.

  It is easy to check that
  \[1_{\cy{u}{\emptyword}}=\left\{\begin{array}{ll}
      \sum_{j=1}^nA(u_{|u|},j)1_{\cyl{j}} & \textrm{if }u\in
      \B(\OSS_A)\\
      0 & \textrm{if } u \notin \B(\OSS_A),
    \end{array} \right. \]
and
  \begin{equation*}
    1_{\cy{u}{v}}=
    \begin{cases}
      1_{\cyl{v}}&\text{if }A(u_1,u_2)=A(u_2,u_3)=\dots
      =A(u_{|u|-1},u_{|u|})=\\ & \hspace{40mm} =A(u_{|u|},v_1)=1,\\
      0&\text{else},
    \end{cases}
  \end{equation*}
  if $v\ne\emptyword$, and that
  \[T_u^*T_u=\left\{\begin{array}{ll}
      \sum_{j=1}^nA(u_{|u|},j)T_jT_j^* & \textrm{if }u\in
      \B(\OSS_A)\\
      0 & \textrm{if } u \notin \B(\OSS_A),
    \end{array} \right. \]
  and
   \begin{equation*}
    T_vT_u^*T_uT_v^*=
    \begin{cases}
      T_vT_v^*&\text{if }A(u_1,u_2)=A(u_2,u_3)=\dots
      =A(u_{|u|-1},u_{|u|})= \\
      & \hspace{40mm}=A(u_{|u|},v_1)=1,\\
      0&\text{else}.
    \end{cases}
  \end{equation*}
  So $\D_{\OSS_A}$ is contained in $\overline{\bigcup_{m\in\N}\D_m}$, and
  $\psi\left(1_{\cy{u}{v}}\right)=T_vT_u^*T_uT_v^*$ for all $u,v\in
  \{1,2,\dotsc ,n\}^*$.

\end{proof}

This result is generalized in \cite{MR2004c:46103}, where it is shown that $\crod$
is isomorphic to a universal Cuntz-Krieger algebra, when $\OSS$ is a sofic shift.

If $A(i,j)=1$ for every $i,j\in\{1,2,\dots,n\}$, then $\Oo_A$, and
hence $\croa$, is the Cuntz algebra $\Oo_n$ which was
originally defined in \cite{MR57:7189}. The Cuntz algebras have
proved to be very important examples in the theory of
$\cs$-algebras, for example in classification of $\cs$-algebras
(see \cite{MR1878882}), and in the study of wavelets (see
\cite{MR1743534}).

\section{Uniqueness and a faithful representation} \label{sec:uniq}

It follows from the universal property of $\crod$ that there exists an
action $\gaug:\T\to\aut(\crod)$ defined by $\gaug_z(S_u)=z^{|u|}S_u$
for every $z\in\T$. This action is known as the \emph{gauge action}.

Let $\F_{\OSS}$ denote the $\cs$-subalgebra of $\crod$ 
generated by $\{S_vS_u^*S_uS_w^*\}_{u,v,w\in\al^*, |v|=|w|}$. It
is not difficult to see that
\begin{multline*}
  \Big\{\sum_{v\in J_-}X_vS_v^*+X_0+\sum_{u\in J_+}S_u
    X_u\mid J_-\text{ and }J_+\text{ are finite subset of }\al^*\\
    \text{
      and }X_0,X_v,X_u\in\F_{\OSS}\text{ for all }v\in J_-,u\in J_+\Big\}
\end{multline*}
is a dense $*$-subalgebra $\crod$. Thus we see that $\F_{\OSS}$ is the
fix point algebra of the gauge action.

If we let
\begin{equation*}
  E(X)=\int_\T\alpha_z(X)dz
\end{equation*}
for every $X\in\crod$, then $E$ is a projection of norm one from
$\crod$ onto $\F_{\OSS}$
satisfying
\begin{eqnarray}
& E(a b c)  =  a E(b)c \mbox{ for all } a,c\in\F_{\OSS}, &
\quad  \mbox{ ({\em module property}) } \\
& E(b^{*}b)  \geq  0, &
\quad   \mbox{ ({\em positivity}) } \\
& E(b^{*}b) = 0 \mbox{ implies that } b=0. & \quad   \mbox{ ({\em
faithfulness}) }
\end{eqnarray}
for all $ b \in \crod$. Thus
\begin{equation*}
  E\Big(\sum_{v\in J_-}X_vS_v^*+X_0+\sum_{u\in J_+}S_uX_u\Big)=X_0
\end{equation*}
for all finite subset $J_-,J_+$ of $\al^*$ and
$X_0,X_v,X_u\in\F_{\OSS},\ v\in J_-,u\in J_+$.

Building on the work done by Matsumoto in \cite{MR98h:46077}, the following Theorem is proved
in \cite{tmciii}:


\begin{theorem}
  Let $\OSS$ be a one-sided shift space, $X$ is a $\cs$-algebra
  generated by partial isometries $\{s_u\}_{u\in\alwords}$, and
  $\phi:\crod\to X$ a $*$-homomorphism such that $\phi(S_u)=s_u$ for
  every $u\in\al^*$. Then the following three statements are
  equivalent:
  \begin{enumerate}
  \item the $*$-homomorphism $\phi:\crod\to X$ is injective,
  \item the restriction of $\phi$ to $\AX$ is injective and there
    exists an action $\gaug:\T\to\aut(X)$ such that $\gaug_z(s_u)=z^{|u|}s_u$
    for every $z\in\T$ and every $u\in\alwords$,
  \item the restriction of $\phi$ to $\AX$ is injective and there
    exists a projection $E$ of norm one from
    $X$ onto $\cs(s_vs_u^*s_us_w^*\mid u,v,w\in\alwords,\ |v|=|w|)$
    satisfying
    \begin{eqnarray*}
      & E(a b c)  =  a E(b)c \mbox{ for all } a,c\in
      \cs(s_vs_u^*s_us_w^*\mid u,v,w\in\alwords,\ |v|=|w|),\\
      & E(b^{*}b)  \geq  0,\\
      & E(b^{*}b) = 0 \mbox{ implies that } b=0, & \quad
    \end{eqnarray*}
    for all $ b \in X$.
  \end{enumerate}
\end{theorem}

As a corollary to this theorem we get:

\begin{corollary} \label{isom}
   Let $\OSS$ be a one-sided shift space. If $X$ is a $\cs$-algebra
  generated by partial isometries $\{s_u\}_{u\in\alwords}$ satisfying:
  \begin{enumerate}
  \item $s_us_v=s_{uv}$ for all $u,v\in\alwords$,
  \item the map
    \begin{equation*}
      1_{\cy{u}{v}}\mapsto s_vs_u^*s_us_v^*,\ u,v\in\alwords
    \end{equation*}
    extends to an injective $*$-homomorphism from $\DX$ to $X$,
  \item there exists an action $\gaug:\T\to\aut(X)$ defined by
    $\gaug_z(s_u)=z^{|u|}s_u$ for every $z\in\T$,
  \end{enumerate}
  then $X$ and $\crod$ are isomorphic by an isomorphism which maps
  $s_u$ to $S_u$ for every $u\in\alwords$.
\end{corollary}

As a consequence of this, we are now able to construct for every
one-sided shift space $\OSS$ a faithful representation of $\crod$
in the following way. Let $\HX$ be the Hilbert space
$l_2(\OSS)\oplus l_2(\Z)$ with orthonormal basis
$(e_x,e_n)_{x\in\OSS,n\in\Z}$, and let for every $u\in\alwords$,
$s_u$ be the operator on $\HX$ defined by:
\begin{equation*}
  s_u(e_x,e_n)=
  \begin{cases}
    (e_{ux},e_{n+|u|})& \text{if } ux\in\OSS,\\
    0& \text{if }ux\notin \OSS.
  \end{cases}
\end{equation*}
It is easy to check that $s_us_v=s_{uv}$ and that
\begin{equation*}
  s_vs_u^*s_us_v^*(e_x,e_n)=
  \begin{cases}
    (e_x,e_n)&\text{if }x\in\cy{u}{v},\\
    0&\text{if }x\notin\cy{u}{v}.
  \end{cases}
\end{equation*}
Thus $\{s_u\}_{u\in\alwords}$ satisfies (1) and (2) of Corollary
\ref{isom}.

If we for every $z\in\T$ let $U_z$ be the operator on $\HX$ defined by
\begin{equation*}
  U_z(e_x,e_n)=z^n(e_x,e_n),
\end{equation*}
then $U_z$ is a unitary operator on $\HX$, and
$U_zs_uU_z^*=z^{|u|}s_u$ for every $u\in\alwords$. Thus
$\{s_u\}_{u\in\alwords}$ also satisfies (3) of Corollary
\ref{isom}, and therefore $S_u\mapsto s_u$ is a faithful
representation of $\crod$.

\begin{definition}
  We say that a one-sided shift space $\OSS$ satisfies condition $(I)$
  if there for every $x\in\OSS$ and every $l\in\No$ exists a
  $y\in\OSS$ such that $\Past_l(x)=\Past_l(y)$ and $x\ne y$.
\end{definition}

One can show that if $\OSS$ satisfies condition (I), then there for
all $\cs$-algebra $X$ generated by partial isometries
$\{s_u\}_{u\in\alwords}$ satisfying:
\begin{enumerate}
\item $s_us_v=s_{uv}$ for all $u,v\in\alwords$,
\item the map
  \begin{equation*}
    1_{\cy{u}{v}}\mapsto s_vs_u^*s_us_v^*,\ u,v\in\alwords
  \end{equation*}
  extends to an injective $*$-homomorphism from $\DX$ to $X$,
\end{enumerate}
exists an action $\gaug:\T\to\aut(X)$ such that
$\gaug_z(s_u)=z^{|u|}s_u$ for every $z\in\T$.

This was first proved by Matsumoto in the case where
$\OSS$ is of the form $\OSS_\Lambda$ for some two-sided shift space
$\Lambda$ in \cite{MR2000d:46082}, where he also discuss several
conditions which are equivalent of condition (I), and this has been
generalized to arbitrary one-sided shift spaces $\OSS$
by the first author in \cite{speciale}.

From this result follows the following theorem:

\begin{theorem}
  Let $\OSS$ be a one-sided shift space which satisfies condition
  (I). If $X$ is a $\cs$-algebra generated by partial isometries
  $\{s_u\}_{u\in\alwords}$ satisfying:
  \begin{enumerate}
  \item $s_us_v=s_{uv}$ for all $u,v\in\alwords$,
  \item the map
    \begin{equation*}
      1_{\cy{u}{v}}\mapsto s_vs_u^*s_us_v^*,\ u,v\in\alwords
    \end{equation*}
    extends to an injective $*$-homomorphism from $\DX$ to $X$,
  \end{enumerate}
  then $X$ and $\crod$ are isomorphic by an isomorphism which maps
  $s_u$ to $S_u$ for every $u\in\alwords$.
\end{theorem}

\section{Properties of $\crod$}
We will in this section shortly describe some of the properties of the
$\cs$-algebra $\crod$.

As mentioned in Remark \ref{corresp}, $\crod$ is isomorphic to the
$\cs$-algebra $\Oo_{\OSS}$ defined in \cite{tmc}, and since $\Oo_{\OSS}$
is the $\cs$-algebra of a separable $\cs$-correspondence over $\DX$
which is separable and commutative and hence nuclear and satisfies the
UCT, the same is the case for the $\cs$-algebra $J_X$ mentioned in
\cite{MR2102572}*{Proposition 8.8}, and thus it follows from
\cite{MR2102572}*{Corollary 7.4 and Proposition 8.8} that $\Oo_{\OSS}$
and hence $\crod$ is nuclear and satisfies the UCT.

\begin{theorem}
  Let $\OSS$ be a one-sided shift space. Then the $\cs$-algebra
  $\crod$ is nuclear and satisfies the UCT.
\end{theorem}

Matsumoto has in \cite{MR2000d:46082} proved the following:
\begin{theorem}
  Let $\Lambda$ be a two-sided shift space. We then have:
  \begin{enumerate}
  \item if $\OSS_\Lambda$ is irreducible in past equivalence, meaning
    that there for every $l\in\No$, every $y\in\OSS_\Lambda$ and every
    sequence $(x_n)_{n\in\N}$ of $\OSS_\Lambda$ such that
    $\Past_n(x_n)=\Past_n(x_{n+1})$ for every $n\in\N$, exist $N\in\N$
    and a $u\in\B(\Lambda)$ such that $\Past_l(y)=\Past_l(ux_{l+N})$,
    then the $\cs$-algebra $\crol$ is simple;
  \item if $\OSS_\Lambda$ is aperiodic in past equivalence, meaning
    that there for any $l\in\No$ exists $N\in\N$ such that for any
    pair $x,y\in\OSS_\Lambda$, exists $u\in\B_N(\Lambda)$ such that
    $\Past_l(y)=\Past_l(ux)$, then the $\cs$-algebra $\crol$ is simple
    and purely infinite.
  \end{enumerate}
\end{theorem}

\section{$\crod$ as an invariant}
We will in this section see that $\crod$ is an invariant for one-sided
conjugacy in the sense that if two one-sided shift spaces $\OSS$ and
$\OSSY$ are conjugate, then $\crod$ and $\croy$ are isomorphic.

This was first proved by Matsumoto in \cite{MR98h:46077} for the special case where
$\OSS=\OSS_{\Lambda}$ and $\OSSY=\OSS_{\Gamma}$ for two two-sided
shift spaces $\Lambda$ and $\Gamma$ satisfying condition (I), and
generalized to the general case in \cite{tmc}. Because of the way we
have constructed $\crod$ in this paper we can very easily prove this
result and even improve it a little bit.

Remember that in $\toepd$, $S^*aS=\tra(a)$ for every $a\in\DX$, so in $\crod$
$\quo(S)^*a\quo(S)=\tra(a)$ for every $a\in\DX$. We will therefore
denote the map
\begin{equation*}
  a\mapsto \quo(S)^*a\quo(S)
\end{equation*}
from $\crod$ to $\crod$ by $\tra$.
We will by $\lambda_{\OSS}$ denote the map
\begin{equation*}
  X\mapsto \left(\sum_{a\in\al}S_a^*\right)X\left(\sum_{b\in\al}S_b\right)
\end{equation*}
from $\F_{\OSS}$ to $\F_{\OSS}$.

\begin{theorem}
  If $\OSS$ and $\OSSY$ are two one-sided shift spaces which are
  conjugate, then there exists a $*$-isomorphism $\Phi$ from $\crod$
  to $\croy$ such that:
  \begin{enumerate}
  \item $\Phi(C(\OSS))=C(\OSSY)$,
  \item $\Phi(\DX)=\D_{\OSSY}$,
  \item $\Phi(\F_{\OSS})=\F_{\OSSY}$,
  \item $\Phi\circ\alpha_{\OSS}=\alpha_{\OSSY}$,
  \item $\Phi\circ\gaug_z=\gaug_z$ for every $z\in\T$,
  \item $\Phi\circ\tra_{\OSS}=\tra_{\OSSY}$,
  \item $\Phi\circ\lambda_{\OSS}=\lambda_{\OSSY}$.
  \end{enumerate}
\end{theorem}

\begin{proof}
  Let $\phi$ be a conjugacy between $\OSSY$ and $\OSS$, and let $\Phi$
  be the map between the bounded functions on $\OSS$ and the bounded functions on
  $\OSSY$ defined by
  \begin{equation*}
    f\mapsto f\circ\phi.
  \end{equation*}
  Then $\Phi(C(\OSS))=C(\OSSY)$, $\Phi\circ\alpha_{\OSS}=\alpha_{\OSSY}\circ\Phi$ and
  $\Phi\circ\tra_{\OSSY}=\tra_{\OSS}\circ\Phi$, and hence
  $\Phi(\D_{\OSSY})=\DX$. Thus it follows from the construction of
  $\crod$ and $\croy$ that there is a $*$-isomorphism from $\crod$ to
  $\croy$ which extends $\Phi$, maps $\quo(S)$ to $\quo(S)$ and
  satisfies $\Phi\circ\alpha_{\OSS}=\alpha_{\OSSY}$. We
  will also denote this $*$-isomorphism by $\Phi$.

  Since the gauge
  action of $\crod$ is characterized by $\gaug_z(f)=f$ for all
  $f\in\DX$ and $\gaug_z(\quo(S))=z\quo(S)$ and the gauge action of $\croy$
  is characterized in the same way, we see that
  $\Phi\circ\gaug_z=\gaug_z$ for every $z\in\T$.

  Since
  $\F_{\OSS}$ is the fix point algebra of the gauge action of $\crod$ and
  $\F_{\OSSY}$ is the fix point algebra of the gauge action of $\croy$,
  we have that $\Phi(\F_{\OSS})=\F_{\OSSY}$.

  Since $\Phi$ maps $\quo(S)$ to $\quo(S)$, we have that
  $\Phi\circ\tra_{\OSS}=\tra_{\OSSY}$.

  Let us denote the function
  \begin{equation*}
     x\mapsto \osh^{-1}\{x\},\ x\in\OSS
  \end{equation*}
  by $f_{\OSS}$ and the function
  \begin{equation*}
    x\mapsto \osh^{-1}\{x\},\ x\in\OSSY
  \end{equation*}
  by $f_{\OSSY}$. We then have that
  \begin{equation*}
    \lambda_{\OSS}(X)=
    \left(\sum_{a\in\al}S_a^*\right)X\left(\sum_{b\in\al}S_b\right)=
    \quo(S)^*\alpha(f_{\OSS})^{1/2}X\alpha(f_{\OSS})^{1/2}\quo(S),
  \end{equation*}
  and since $\Phi(f_{\OSS})=f_{\OSSY}$, we have that
  $\Phi\circ\lambda_{\OSS}=\lambda_{\OSSY}$.
\end{proof}

If two two-sided shift spaces $\Lambda$ and $\Gamma$ are flow equivalent,
then the corresponding one-sided shift spaces $\OSS_{\Lambda}$ and
$\OSS_{\Gamma}$ are not necessarily conjugate, so we cannot expect
that $\crol$ and $\cro{\D_{\OSS_{\Gamma}}\mspace{-5mu}}{\alpha}{\tra}$ are isomorphic
(and there are examples of two two-sided shift spaces $\Lambda$ and
$\Gamma$, such that $\Lambda$ and $\Gamma$ are conjugate and hence
flow equivalent, but $\crol$
and $\cro{\D_{\OSS_{\Gamma}}\mspace{-5mu}}{\alpha}{\tra}$ are not isomorphic), but it
turns out that $\crol\otimes\K$ and
$\cro{\D_{\OSS_{\Gamma}}\mspace{-5mu}}{\alpha}{\tra}\otimes\K,$ where $\K$ is the $C^*$-algebra of compact operators on a separable Hilbert space. 
This has been proved by Matsumoto in \cite{MR2001e:46115}
for $\Lambda$ and $\Gamma$ satisfying condition (I), and in generality
in \cite{tmciii}.

\section{The $K$-theory of $\crod$}
Since $K_0$ and $K_1$ are invariants of a $\cs$-algebra, it follows
from the previous section that $K_0(\crod)$, $K_1(\crod)$ and
$K_0(\F_{\OSS})$ are invariants of $\OSS$. We will in this section
present formulas based on $l$-past equivalence for these
invariants. This was done in
\cites{MR2002h:19004,MR2000d:46082,MR2000e:46087} for the case of
one-sided shift spaces of the form $\OSS_\Lambda$, where $\Lambda$ is
a two-sided shift space and generalized to the general case in \cite{speciale}.

One can directly from these formulas prove that there are
invariants of $\OSS$ without involving $\cs$-algebras. This is done
(for one-sided shift spaces of the form $\OSS_\Lambda$, where
$\Lambda$ is a two-sided shift space) in
Matsumoto's outstanding paper \cite{MR2000h:37013}, where also other
invariants of shift spaces are presented.

Let $\OSS$ be a one-sided shift space. We let for each $l\in\No$,
$m(l)$ be the number of $l$-past equivalence classes and we denote the
$l$-past equivalence classes by
$\E1^l,\E_2^l,\dotsc,\E_{m(l)}^l$.
For each  $l\in \No,\ j\in \{1,2,\dots ,m(l)\}$ and $i\in \{1,2,\dots
,m(l+1)\}$ we let
\[ I_l(i,j)=\left\{ \begin{array}{ll}
        1 & \textrm{if }
        \E_i^{l+1}\subseteq \E_j^l \\
        0 & \textrm{else.}
\end{array} \right.\]

Let $F$ be a finite set and $i_0\in F$. Then we denote by $e_{i_0}$
the element in $\Z^F$ for which
\[e_{i_0}(i)=\left\{ \begin{array}{ll} 1 & \textrm{if } i=i_0\\
0 & \textrm{else.}
\end{array} \right. \]
For $0\le k \le l$ let $M_k^l=\{i\in\{1,2\dots ,m(l)\}\mid
\Past_k(\E_i^l)\ne \emptyset\}$. Since $i\in M_k^{l+1}$ if
$j\in M_k^l$ and $I_l(i,j)=1$, there exists a
positive, linear map from $\Z^{M_k^l}$ to $\Z^{M_k^{l+1}}$ given by
\[e_j\mapsto \sum_{i\in M_k^{l+1}}I_l(i,j)e_i.\]
We denotes this map by $I_k^l$.



For a subset $\E$ of $\OSS$ and a $u\in \alwords$ we let
$u\E=\{u x\in \OSS \mid x\in \E\}$.
For each $l\in \No,\ j\in \{1,2,\dots ,m(l)\},\ i\in \{1,2,\dots
,m(l+1)\}$ and $a\in \al$ we let
\[ A_l(i,j,a)=\left\{ \begin{array}{ll}
        1 & \textrm{if }
        \emptyset \ne a\E_i^{l+1}\subseteq \E_j^l \\
        0 & \textrm{else.}
\end{array} \right.\]

Let $0\le k\le l$. If $j\in M_k^l$ and there exists an $a\in \al$ such that $A_l(i,j,a)=1$,
then $i\in M_{k+1}^{l+1}$. Hence
there exists a positive,  linear map from $\Z^{M_k^l}$ to
$\Z^{M_{k+1}^{l+1}}$ given by
\[e_j\mapsto \sum_{i\in M_{k+1}^{l+1}}\sum_{a\in \al}A_l(i,j,a)e_i.\]
We denote this map by $A_k^l$.

Then we have:

\begin{lemma} \label{kisomt}
Let $0\le k \le l$. Then the following
diagram commutes:
\begin{displaymath}
\xymatrix{ \Z^{M_k^l} \ar[r]^{I_k^l} \ar[d]_{A_k^l} & \Z^{M_k^{l+1}} \ar[d]^{A_k^{l+1}} \\
        \Z^{M_{k+1}^{l+1}} \ar[r]^{I_{k+1}^{l+1}} &
        \Z^{M_{k+1}^{l+2}}. }
\end{displaymath}
\end{lemma}

\begin{proof}
  Let $j\in M_k^l$, $h\in M_{k+1}^{l+2}$ and $a\in \al$. If $\emptyset
  \ne a\E_h^{l+2}\subseteq \E_j^l$, then there exists exactly one
  $i\in M_k^{l+1}$ such that $\E_i^{l+1}\subseteq \E_j^l$ and
  $\emptyset \ne a\E_h^{l+2}\subseteq \E_i^{l+1}$; and there exists
  exactly one $i'\in M_{k+1}^{l+1}$ such that $\E_h^{l+2}\subseteq
  \E_{i'}^{l+1}$ and $\emptyset \ne a\E_{i'}^{l+1}\subseteq \E_j^l$;
  and if $a\E_h^{l+2}=\emptyset $ or $a\E_h^{l+2}\nsubseteq \E_j^l$
  then there does not exists a $i\in M_k^{l+1}$ such that
  $\E_i^{l+1}\subseteq \E_j^l$ and $\emptyset \ne a\E_h^{l+2}\subseteq
  \E_i^{l+1}$; and there does not exists a $i'\in M_{k+1}^{l+1}$ such
  that $\E_h^{l+2}\subseteq \E_{i'}^{l+1}$ and $\emptyset \ne
  a\E_{i'}^{l+1}\subseteq \E_j^l$. Hence
  \[\sum_{i\in M_k^{l+1}}A_{l+1}(h,i,a)I_l(i,j)=\sum_{i\in
    M_{k+1}^{l+1}}I_{l+1}(h,i)A_l(i,j,a) .\]
  So
  \begin{eqnarray*}
    A_k^{l+1}(I_k^l(e_j))&=& A_k^{l+1}\left(\sum_{i\in
        M_k^{l+1}}I_l(i,j)e_i\right)\\
    &=& \sum_{h\in M_{k+1}^{l+2}}\sum_{a\in
      \al}A_{l+1}(h,i,a)\sum_{i\in M_k^{l+1}}I_l(i,j)e_h\\
    &=& \sum_{h\in M_{k+1}^{l+2}}\sum_{i\in M_{k+1}^{l+1}}\sum_{a\in
      \al}I_{l+1}(h,i)A_l(i,j,a)e_h\\
    &=& I_{k+1}^{l+1}\left(\sum_{i\in M_{k+1}^{l+1}}\sum_{a\in
        \al}A_l(i,j,a)e_i\right) \\
    &=& I_{k+1}^{l+1}(A_k^l(e_j))
  \end{eqnarray*}
  for every $j\in M_k^l$.
  Hence the diagram commutes.
\end{proof}
For $k\in \No$ the inductive limit
$\limm(\Z^{M_k^l},(\Z^+)^{M_k^l},I_k^l)$ will be denoted by
$(\Z_{\OSS_k},\Z_{\OSS_k}^+)$. It follows from Lemma \ref{kisomt}
that the $A_k^l$'s induce a positive, linear map $A_k$ from
$\Z_{\OSS_k}$ to $\Z_{\OSS_{k+1}}$.

Let $0\le k<l$. Denote by $\delta_k^l$ the linear map from $\Z^{M_k^l}$ to
$\Z^{M_{k+1}^l}$ given by
\begin{equation*}
  e_j\mapsto
  \begin{cases}
    e_j&\text{if }j\in M_{k+1}^l,\\
    0&\text{if }j\notin M_{k+1}^l,
  \end{cases}
\end{equation*}
for $j\in M_k^l$. It is easy to check that the following diagram
commutes
\begin{displaymath}
\xymatrix{ \Z^{M_k^l}
  \ar[r]^{\delta_k^l} \ar[d]_{I_k^l} & \Z^{M_{k+1}^l} \ar[d]^{I_{k+1}^l} \\
\Z^{M_k^{l+1}} \ar[r]^{\delta_k^{l+1}} & \Z^{M_{k+1}^{l+1}}. }
\end{displaymath}
Thus the $\delta_k^l$'s induce a positive, linear map from
$\Z_{\OSS_k}$ to $\Z_{\OSS_{k+1}}$ which we denote by
$\delta_k$. Since the diagram
\begin{displaymath}
\xymatrix{ \Z^{M_k^l}
  \ar[r]^{\delta_k^l} \ar[d]_{A_k^l} & \Z^{M_{k+1}^l} \ar[d]^{A_{k+1}^l} \\
\Z^{M_{k+1}^{l+1}} \ar[r]^{\delta_{k+1}^{l+1}} & \Z^{M_{k+2}^{l+1}}  }
\end{displaymath}
commutes for every $0\le k<l$, the diagram
\begin{displaymath}
\xymatrix{ \Z_{\OSS_k} \ar[r]^{\delta_k} \ar[d]_{A_k} &
  \Z_{\OSS_{k+1}} \ar[d]^{A_{k+1}} \\
\Z_{\OSS_{k+1}} \ar[r]^{\delta_{k+1}} & \Z_{\OSS_{k+2}}}
\end{displaymath}
commutes.

We denote the inductive
limit $\limm(\Z_{\OSS_k},\Z_{\OSS_k}^+,A_k)$ by
$(\Delta_{\OSS},\Delta_{\OSS}^+)$. Since the previous diagram
commutes, the $\delta_k$'s induce a positive, linear map from
$\Delta_{\OSS}$ to $\Delta_{\OSS}$ which we denote by
$\delta_{\OSS}$.

\begin{theorem}
  For every one-sided shift space $\OSS$ is
  \[(K_0(\F_{\OSS}),K_0^+(\F_{\OSS}),(\lambda_{\OSS})_*)\cong
  (\Delta_{\OSS},\Delta_{\OSS}^+,\delta_{\OSS}),\]
 or more precisely, 
  the map $[S_u 1_{\E^l_i}S_v^*]_0\mapsto e_i\in \Z^{M_k^l}$
  extends to an isomorphism from
  $(K_0(\F_{\OSS}),K_0^+(\F_{\OSS}),(\lambda_{\OSS})_*)$ to
  $(\Delta_{\OSS},\Delta_{\OSS}^+,\delta_{\OSS})$.
\end{theorem}

Denote for every $l\in \No$ by $B^l$ the linear map from $\Z^{M_1^l}$
to $\Z^{m(l+1)}$ given by
\[e_j\mapsto \sum_{i=1}^{m(l+1)}\left(I_l(i,j)-\sum_{a\in \al}A_l(i,j,a)\right)e_i.\]

\noindent One can easily check that the following diagram commutes for every
$l\in \No$.
\begin{displaymath}
\xymatrix{ \Z^{M_1^l}
  \ar[r]^{B^l} \ar[d]_{I_1^l} & \Z^{m(l+1)} \ar[d]^{I_0^{l+1}} \\
\Z^{M_1^{l+1}} \ar[r]^{B^{l+1}} & \Z^{m(l+2)}. }
\end{displaymath}
Hence the $B^l$'s induce a linear map $B$ from $\Z_{\Lambda_1}$ to
$\Z_{\Lambda_0}$.



\begin{theorem} \label{heks}
Let $\Lambda$ be a one-sided shift space. Then
\[K_0(\Oo_\Lambda)\cong \Z_{\Lambda_0}/B\Z_{\Lambda_1},\]
and
\[K_1(\Oo_\Lambda)\cong \ker(B).\]
More precisely: The map
\[[1_{\E_i^l}]_0\mapsto e_i\in \Z^{m(l)}\]
induces an isomorphism from $K_0(\Oo_\Lambda)$ to
$\Z_{\Lambda_0}/B\Z_{\Lambda_1}$.
\end{theorem}


\section{The ideal structure of $\crod$}
We will in this section describe the structure of the gauge invariant
ideals of $\crod$. By an ideal we will in this paper always mean a
closed two-sided ideal, and by a gauge invariant ideal, we mean an
ideal $I$ such that $\gaug_z(I)\subseteq I$ for every $z\in\T$.

The lattice of the gauge invariant ideals of $\crod$ has been
described by Matsumoto in \cite{MR2000d:46082} in the case where
$\OSS$ is of the form $\OSS_\Lambda$ for some two-sided shift space
$\Lambda$ and this has been generalized to arbitrary
one-sided shift spaces $\OSS$ by the first author in
\cite{speciale}. We will here reformulate the description a bit.
\begin{theorem}
  Let $\OSS$ be a one-sided shift space. Then there exist between each
  pair of the following lattices an ordering preserving bijective map:
  \begin{enumerate}
  \item the lattice of gauge invariant ideals of $\crod$,
  \item the lattice of ideals $J$ of $\F_{\OSS}$, such that $S_uXS_u^*,
    S_u^*XS_u\in J$ for every $u\in\alwords$ and every $X\in J$,
  \item the lattice of ideals $I$ of $\AX$, such that $S_u^*XS_u\in I$
    for every $u\in\alwords$ and every $X\in I$,
  \item the lattice of order ideals of $\Delta_{\OSS}$ invariant under
    $\delta_{\OSS}$,
  \item the lattice of subset $A$ of $\OSS$, such that
    $\osh(A)\subseteq A$ and $\forall x\in A \ \exists
    l\in\No:\Past_l(x)\subseteq A$.
  \end{enumerate}
\end{theorem}

\section{Examples}

If $\Lambda$ is a two-sided shift space, then as explained before we
can associate to it the $\cs$-algebra $\crol$, but we of course also
look at the $\cs$-crossed product $C(\Lambda)\rtimes_\phi \Z$,
where $\phi:C(\Lambda)\to C(\Lambda)$ is the map
\begin{equation*}
  f\mapsto f\circ\tsh.
\end{equation*}
It is proved in \cite{tmcii} that if $\Lambda$ satisfy the condition
\begin{equation*}
  (*): \text{There exists for every }u\in\lanl \text{ an
  }x\in\OSS_\Lambda \text{ such that }\Past_{|u|}(x)=\{u\},
\end{equation*}
then $C(\Lambda)\rtimes_\phi \Z$ is a quotient of $\crol$. This is
used in \cite{tmcseii} and \cite{tmcseiv} to relate the $K$-theory of
$\crol$ to the $K$-theory of $C(\Lambda)\rtimes_\phi \Z$ for these
shift spaces, and in \cite{MR2085388} to present $K_0(\crol)$, for a
two-sided shift space $\Lambda$ associated to an aperiodic and
primitive substitution, as a stationary inductive limit of a system
associated to an integer matrix defined from combinatorial data which
can be computed in an algorithmic way (cf.\ \cite{tmcsei} and
\cite{tmcse:program}).

In \cite{MR2088934}, Matsumoto has taken a closer look at $\crod$ in
the case where $\OSS$ is the Motzkin shift, and in \cite{MR1716953} he
examines $\crod$ for the context-free shift. In \cite{MR2029162}
$\crod$ is examined for the Dyck shift, and in \cite{MR1645334} $\crod$
is examined for a class of shift spaces called $\beta$-shifts.

{\bf Acknowledgement. }
{\footnotesize We are grateful to Johan {\"O}inert for many useful comments.}

\end{document}